\documentclass[final,3p,times]{elsarticle}

\newtheorem{theorem}{Theorem}[section]
\newtheorem{lemma}[theorem]{Lemma}

\newtheorem{proposition}[theorem]{Proposition}

\newdefinition{definition}[theorem]{Definition}
\newdefinition{example}[theorem]{Example}
\newdefinition{xca}[theorem]{Exercise}

\newdefinition{remark}[theorem]{Remark}
\newproof{proof}{Proof}
\newproof{vrfctn}{Verification}

\usepackage{times}
\usepackage{amssymb}
\usepackage{amsfonts}
\usepackage{amsmath}
\usepackage{euscript}
\usepackage{color}
\usepackage{graphicx}
\usepackage{rotating}
\usepackage{bm}
\usepackage{curves}
\usepackage[svgnames,dvipsnames]{pstricks}
\usepackage{lettrine}
\usepackage{mathrsfs}
\usepackage{booktabs}
\DeclareGraphicsExtensions{{},.eps}
\providecommand{\tsc}[1]{{\text{\sc#1}}}

\providecommand{\etal}{et al.}

\providecommand{\lemref}[1]{{\textup{(Lemma~\ref{#1})}}}

\providecommand{\prpref}[1]{{\textup{(Proposition~\ref{#1})}}}
\providecommand{\defref}[1]{{\textup{(Definition~\ref{#1})}}}

\providecommand{\rmkref}[1]{{\textup{(Remark~\ref{#1})}}}

\providecommand{\lemrefnp}[1]{{\textup{Lemma~\ref{#1}}}}

\providecommand{\prprefnp}[1]{{\textup{Proposition~\ref{#1}}}}
\providecommand{\defrefnp}[1]{{\textup{Definition~\ref{#1}}}}

\providecommand{\rmkrefnp}[1]{{\textup{Remark~\ref{#1}}}}

\providecommand{\eqrefsatob} [2]{\textup{(\ref{#1}--\ref{#2})}}
\providecommand{\eqrefsab}   [2]{\textup{(\ref{#1}, \ref{#2})}}
\providecommand{\eqrefsabc}  [3]{\textup{(\ref{#1}, \ref{#2}, \ref{#3})}}

\providecommand{\parrefnp}[1]{{\textup{\S\ref{#1}}}}
\providecommand{\const}{{\rm const}}

\begin{document}

\begin{frontmatter}

\title{A general recurrence relation for the weight-functions in\\
       M\"uhlbach-Neville-Aitken representations\\
       with application to WENO interpolation and differentiation}

\author{G.A. Gerolymos}
\ead{georges.gerolymos@upmc.fr}
\address{Universit\'e Pierre-et-Marie-Curie \textup{(}\tsc{upmc}\textup{)}, 4 place Jussieu, 75005 Paris, France}
\journal{Appl. Math. Comp.}


\begin{abstract}
In several applications, such as \tsc{weno} interpolation and reconstruction [Shu C.W.: {\em SIAM Rev.} {\bf 51} (2009) 82--126],
we are interested in the analytical expression of the weight-functions which allow the representation of the approximating function
on a given stencil (Chebyshev-system) as the weighted combination of the corresponding approximating functions on substencils (Chebyshev-subsystems).
We show that the weight-functions in such representations [M\"uhlbach G.: {\em Num. Math.} {\bf 31} (1978) 97--110] can be generated by a general
recurrence relation based on the existence of a $1$-level subdivision rule. As an example of application we apply this recurrence to the
computation of the weight-functions for Lagrange interpolation [Carlini E., Ferretti R., Russo G.: {\em SIAM J. Sci. Comp.} {\bf 27} (2005) 1071--1091]
for a general subdivision of the stencil $\{x_{i-M_-},\cdots,x_{i+M_+}\}$ of $M+1:=M_-+M_++1$ distinct ordered points
into $K_\mathrm{s}+1\leq M:=M_-+M_+>1$ (Neville) substencils $\{x_{i-M_-+k_\mathrm{s}},\cdots,x_{i+M_+-K_\mathrm{s}+k_\mathrm{s}}\}$ ($k_\mathrm{s}\in\{0,\cdots,K_\mathrm{s}\}$)
all containing the same number of $M-K_\mathrm{s}+1$ points but each shifted by 1 cell with respect to its neighbour,
and give a general proof for the conditions of positivity of the weight-functions (implying convexity of the combination),
extending previous results obtained for particular stencils and subdvisions [Liu Y.Y., Shu C.W., Zhang M.P.: {\em Acta Math. Appl. Sinica} {\bf 25} (2009) 503--538].
Finally, we apply the recurrence relation to the representation by combination of substencils of derivatives of arbitrary order of the Lagrange interpolating polynomial.
\end{abstract}

\end{frontmatter}

%
%
%
%
%
%
%
%
%
\section{Introduction}\label{GRRWFsMNARsAWENOID_s_I}
%
%
%
%
%
%
%
%
%

The Neville-Aitken algorithm \cite[pp. 204--209]{Henrici_1964a} constructs the interpolating polynomial on $\{x_{i-M_-},\cdots,x_{i+M_+}\}$, by recursive combination of the
interpolating polynomials on substencils, with weight-functions which are also polynomials of $x$ \cite[pp. 204--209]{Henrici_1964a}.
Carlini~\etal~\cite{Carlini_Ferretti_Russo_2005a}, working on the Lagrange interpolating polynomial in the context of centered (central) \tsc{weno} schemes \cite{Shu_2009a},
recognized the connexion between the Neville algorithm \cite[pp. 207--208]{Henrici_1964a} and the determination of the optimal \cite{Shu_2009a} weight-functions.
%
\begin{definition}[{\rm Stencil}]
\label{Def_GRRWFsMNARsAWENOID_s_I_001}
Let
\begin{equation}
\tsc{x}_{i-M_-,i+M_+}:=\{x_{i-M_-},\cdots,x_{i+M_+}\}\subset\mathbb{R} \qquad\left\{\begin{array}{l}M_\pm\in\mathbb{Z}:M:=M_-+M_+\geq0\\
                                                                                                    x_{i-M_-}<\cdots<x_{i+M_+}\;\forall M>0\\\end{array}\right.
                                                                                                                                    \label{Eq_Def_GRRWFsMNARsAWENOID_s_I_001_001}
\end{equation}
be a set of $M+1$ distinct ordered real points.
\qed
\end{definition}
%
%
\begin{definition}[{\rm Neville substencils}]
\label{Def_GRRWFsMNARsAWENOID_s_I_002}
\begin{subequations}
                                                                                                                                    \label{Eq_Def_GRRWFsMNARsAWENOID_s_I_002_001}
Let $\tsc{x}_{i-M_-,i+M_+}$ be a stencil \defref{Def_GRRWFsMNARsAWENOID_s_I_001} and assume $M\geq2$ in \eqref{Eq_Def_GRRWFsMNARsAWENOID_s_I_001_001}.
The $K_{\rm s}+1>1$ substencils
\begin{equation}
\tsc{x}_{i-M_-+k_\mathrm{s},i+M_+-K_\mathrm{s}+k_\mathrm{s}}:=\left\{x_{i-M_-+k_\mathrm{s}},\cdots,x_{i+M_+-K_\mathrm{s}+k_\mathrm{s}}\right\}
\qquad\left\{\begin{array}{l}{\mathbb N}\ni K_{\rm s}\leq M-1:=M_-+M_+-1\geq1\\
                             k_\mathrm{s}\in\{0,\cdots,K_{\rm s}\}           \\\end{array}\right.
                                                                                                                                    \label{Eq_Def_GRRWFsMNARsAWENOID_s_I_002_001a}
\end{equation}
each of which contains $M-K_\mathrm{s}+1\geq2$ distinct ordered points and which satisfy
\begin{alignat}{6}
&
\bigcup_{k_\mathrm{s}=0}^{K_\mathrm{s}}\tsc{x}_{i-M_-+k_{\rm s},i+M_+-K_\mathrm{s}+k_\mathrm{s}}=\tsc{x}_{i-M_-,i+M_+}
&
                                                                                                                                    \label{Eq_Def_GRRWFsMNARsAWENOID_s_I_002_001b}\\
&
\tsc{x}_{i-M_-+k_\mathrm{s}+1,i+M_+-K_\mathrm{s}+k_\mathrm{s}+1}=\Big(\tsc{x}_{i-M_-+k_\mathrm{s},i+M_+-K_\mathrm{s}+k_\mathrm{s}}\setminus\{x_{i-M_-+k_\mathrm{s}}\}\Big)\cup\{x_{i+M_+-K_\mathrm{s}+k_\mathrm{s}+1}\}
\quad\forall k_\mathrm{s}\in\{0,\cdots,K_{\rm s}-1\}
&
                                                                                                                                    \label{Eq_Def_GRRWFsMNARsAWENOID_s_I_002_001c}
\end{alignat}
correspond to the $K_\mathrm{s}$-level subdivision of $\tsc{x}_{i-M_-,i+M_+}$ to substencils of equal length,
each obtained from its left neighbour by deleting the leftmost point and adding 1 point to the right \eqref{Eq_Def_GRRWFsMNARsAWENOID_s_I_002_001c}.
\end{subequations}
\qed
\end{definition}
%
The optimal weight-functions $\sigma_{I,\tsc{x}_{i-M_-,i+M_+},K_\mathrm{s},k_\mathrm{s}}(x)$ in \tsc{weno} interpolation satisfy \cite{Carlini_Ferretti_Russo_2005a,
                                                                                                                                       Shu_2009a,
                                                                                                                                       Liu_Shu_Zhang_2009a}
\begin{subequations}
                                                                                                                                    \label{Eq_GRRWFsMNARsAWENOID_s_I_001}
\begin{alignat}{6}
p_{I,\tsc{x}_{i-M_-,i+M_+}}(x;f)= \sum_{k_\mathrm{s}=0}^{K_\mathrm{s}}\sigma_{I,\tsc{x}_{i-M_-,i+M_+},K_\mathrm{s},k_\mathrm{s}}(x)
                                  \;p_{I,\tsc{x}_{i-M_-+k_\mathrm{s},i+M_+-K_\mathrm{s}+k_\mathrm{s}}}(x;f)                            &\qquad \forall x\in\mathbb{R}
                                                                                                                                    \label{Eq_GRRWFsMNARsAWENOID_s_I_001a}\\
                                  \sum_{k_\mathrm{s}=0}^{K_\mathrm{s}}\sigma_{I,\tsc{x}_{i-M_-,i+M_+},K_\mathrm{s},k_\mathrm{s}}(x) = 1&\qquad \forall x\in\mathbb{R}
                                                                                                                                    \label{Eq_GRRWFsMNARsAWENOID_s_I_001b}
\end{alignat}
\end{subequations}
where $p_{I,\tsc{x}_{i-M_-,i+M_+}}(x;f)$ is the Lagrange interpolating polynomial \cite[pp. 183--189]{Henrici_1964a}
of the real function $f:\mathbb{R}\to\mathbb{R}$ on the stencil $\tsc{x}_{i-M_-,i+M_+}$ \defref{Def_GRRWFsMNARsAWENOID_s_I_001}.
The optimal weight-functions $\sigma_{I,\tsc{x}_{i-M_-,i+M_+},K_\mathrm{s},k_\mathrm{s}}(x)$ \eqref{Eq_GRRWFsMNARsAWENOID_s_I_001} correspond to the weight-functions in M\"uhlbach's theorem \cite[Theorem 2.1, p. 100]{Muhlbach_1978a},
where they were expressed in terms of quotients of determinants of interpolation-error functions, directly obtained by the Cramer solution \cite[Proposition 5.1.1, p. 72]{Allaire_Kaber_2008a}
of error-eliminating linear systems \cite[(13), p. 8489]{Gerolymos_Senechal_Vallet_2009a}.
Since \tsc{weno} approaches are based on nonlinear, with respect to the function $f(x)$, modifications of the optimal weight-functions \eqref{Eq_GRRWFsMNARsAWENOID_s_I_001},
we are particularly interested in analytical explicit expressions of the weight-functions.

Carlini~\etal~\cite[(3.6,4.10), pp. 1074--1079]{Carlini_Ferretti_Russo_2005a}
gave the expression of the polynomial weight-functions $\sigma_{I,\tsc{x}_{i-(r-1),i+r},r-1,k_\mathrm{s}}(x)$ for the $(K_\mathrm{s}=r-1)$-level subdivision \defref{Def_GRRWFsMNARsAWENOID_s_I_002}
of $\tsc{x}_{i-(r-1),i+r}$ \defref{Def_GRRWFsMNARsAWENOID_s_I_001}.
This result was also confirmed by Liu~\etal~\cite[(2.2), p. 506]{Liu_Shu_Zhang_2009a} who further gave the analytical expression \cite[(2.18), p. 511]{Liu_Shu_Zhang_2009a}
for the polynomial weight-functions $\sigma_{I,\tsc{x}_{i-r,i+r},r,k_\mathrm{s}}(x)$ for the $(K_\mathrm{s}=r)$-level subdivision \defref{Def_GRRWFsMNARsAWENOID_s_I_002}
of $\tsc{x}_{i-r,i+r}$ \defref{Def_GRRWFsMNARsAWENOID_s_I_001}.
In both cases it was observed \cite{Carlini_Ferretti_Russo_2005a,
                                    Liu_Shu_Zhang_2009a}
that $\forall x\in[x_{i-1},x_{i+1}]$ the linear weight-functions are positive ($\geq0$), so that, because of the consistency relation \eqref{Eq_GRRWFsMNARsAWENOID_s_I_001b},
the combination \eqref{Eq_GRRWFsMNARsAWENOID_s_I_001a} of substencils is convex $\forall x\in[x_{i-1},x_{i+1}]$.

The purpose of the present note is to prove \lemref{Lem_GRRWFsMNARsAWENOID_s_GRRWFs_001} a general recurrence relation \cite[(70), p. 299]{Gerolymos_2011a} for weight-functions of an arbitrary family of functions,
for which the $(K_\mathrm{s}=1)$-level subdivision \defref{Def_GRRWFsMNARsAWENOID_s_I_002} is possible. As an example of application we apply this relation to the Lagrange interpolating
polynomial \prpref{Prp_GRRWFsMNARsAWENOID_s_ALIP_001}, for an arbitrary level of subdivision \defref{Def_GRRWFsMNARsAWENOID_s_I_002} of a general stencil \defref{Def_GRRWFsMNARsAWENOID_s_I_001}.
The explicit expression of the weight-funcions developed in \prprefnp{Prp_GRRWFsMNARsAWENOID_s_ALIP_001} is used \prpref{Prp_GRRWFsMNARsAWENOID_s_ALIP_002} to study the convexity
of representation \eqref{Eq_GRRWFsMNARsAWENOID_s_I_001}. Then, we apply the general recurrence relation \lemref{Lem_GRRWFsMNARsAWENOID_s_GRRWFs_001} to determine the weight-functions
for the representation of the $n$-derivative of the Lagrange interpolating polynomial by combination of substencils \prpref{Prp_GRRWFsMNARsAWENOID_s_AnDLIP_001}.

%
%
%
%
%
%
%
%
%
\section{General recurrence relation for weight-functions}\label{GRRWFsMNARsAWENOID_s_GRRWFs}
%
%
%
%
%
%
%
%
%

We start by considering a more abstract case, where a general family of functions $p_{M_-,M_+}(x)$ depending on 2 integer indices $M_\pm\in\mathbb{Z}:M_-+M_+\geq1$
(which in practical applications may correspond to stencils; \defrefnp{Def_GRRWFsMNARsAWENOID_s_I_001}),
are equipped with a $1$-level subdivion property, and develop a general recurrence relation for the weight-functions.
%
\begin{lemma}[{\rm Recursive generation of weight-functions}]
\label{Lem_GRRWFsMNARsAWENOID_s_GRRWFs_001}
\begin{subequations}
                                                                                                                                    \label{Eq_Lem_GRRWFsMNARsAWENOID_s_GRRWFs_001_001}
Let $p_{M_-,M_+}(x)$ be a family of real functions
\begin{alignat}{6}
p_{M_-,M_+}:\mathbb{R}\to\mathbb{R}\qquad\forall M_\pm\in\mathbb{Z}\;:\;M:=M_-+M_+\geq1\Longrightarrow M_+>-M_-
                                                                                                                                    \label{Eq_Lem_GRRWFsMNARsAWENOID_s_GRRWFs_001_001a}
\end{alignat}
and assume that there exists an associated family of real weight-functions $\sigma_{M_-,M_+,1,0}(x)$ \textup{(}also defining $\sigma_{M_-,M_+,1,1}(x):=1-\sigma_{M_-,M_+,1,0}(x)$\textup{)}
\begin{alignat}{6}
&\sigma_{M_-,M_+,1,k_\mathrm{s}}:\mathbb{R}\setminus\mathcal{S}_{\sigma_{M_-,M_+,1}}\to\mathbb{R}\qquad\left\{\begin{array}{l}\forall M_\pm\in\mathbb{Z}\;:\;M:=M_-+M_+\geq2\\
                                                                                                                               \forall k_\mathrm{s}\in\{0,1\}\\\end{array}\right.
                                                                                                                                    \label{Eq_Lem_GRRWFsMNARsAWENOID_s_GRRWFs_001_001b}\\
&\sigma_{M_-,M_+,1,0}(x)+\sigma_{M_-,M_+,1,1}(x)=1\qquad\left\{\begin{array}{l}\forall x\in\mathbb{R}\setminus\mathcal{S}_{\sigma_{M_-,M_+,1}}\\
                                                                               \forall M_\pm\in\mathbb{Z}\;:\;M:=M_-+M_+\geq2\\\end{array}\right.
                                                                                                                                    \label{Eq_Lem_GRRWFsMNARsAWENOID_s_GRRWFs_001_001c}
\end{alignat}
defined everywhere in $\mathbb{R}$ except for a finite set of isolated points $\mathcal{S}_{\sigma_{M_-,M_+,1}}\subset\mathbb{R}$, which may be empty, such that
\begin{alignat}{6}
p_{M_-,M_+}(x)=\sigma_{M_-,M_+,1,0}(x)\;p_{M_-,M_+-1}(x)+\sigma_{M_-,M_+,1,1}(x)\;p_{M_--1,M_+}(x)\qquad\left\{\begin{array}{l}\forall x\in\mathbb{R}\setminus\mathcal{S}_{\sigma_{M_-,M_+,1}}\\
                                                                                                                               \forall M_\pm\in\mathbb{Z}\;:\;M:=M_-+M_+\geq2\\\end{array}\right.
                                                                                                                                    \label{Eq_Lem_GRRWFsMNARsAWENOID_s_GRRWFs_001_001d}
\end{alignat}
Then the family of weight-functions defined recursively by
\begin{alignat}{6}
\sigma_{M_-,M_+,K_\mathrm{s},k_\mathrm{s}}(x):= \sum_{\ell_\mathrm{s}=\max(0,k_\mathrm{s}-1)}^{\min(K_\mathrm{s}-1,k_\mathrm{s})}
                                                \sigma_{M_-                ,M_+                                 ,K_\mathrm{s}-1,            \ell_\mathrm{s}}(x)\;
                                                \sigma_{M_--\ell_\mathrm{s},M_+-(K_\mathrm{s}-1)+\ell_\mathrm{s},1             ,k_\mathrm{s}-\ell_\mathrm{s}}(x)
\qquad\left\{\begin{array}{l}\forall M_\pm\in\mathbb{Z}\;:\;M_-+M_+\geq2     \\
                             \forall K_\mathrm{s}\leq M_-+M_+-1              \\
                             \forall k_\mathrm{s}\in\{0,\cdots,K_\mathrm{s}\}\\\end{array}\right.
                                                                                                                                    \label{Eq_Lem_GRRWFsMNARsAWENOID_s_GRRWFs_001_001e}
\end{alignat}
satisfies
\begin{alignat}{6}
&p_{M_-,M_+}(x)=\sum_{k_\mathrm{s}=0}^{K_\mathrm{s}}\sigma_{M_-,M_+,K_\mathrm{s},k_\mathrm{s}}(x)\;p_{M_--k_\mathrm{s},M_+-K_\mathrm{s}+k_\mathrm{s}}(x)&
\quad\left\{\begin{array}{l}\displaystyle\forall x\in{\mathbb R}\setminus\bigcup_{   L_\mathrm{s}=0}^{K_\mathrm{s}-1}\displaylimits
                                                                         \bigcup_{\ell_\mathrm{s}=0}^{L_\mathrm{s}  }\mathcal{S}_{\sigma_{M_--\ell_\mathrm{s},M_+-L_\mathrm{s}+\ell_\mathrm{s},1}}\\
                                                                                                                                                                                                  \\
                                          \forall M_\pm\in\mathbb{Z}\;:\;M:=M_-+M_+\geq2                                                                                                          \\
                                          \forall K_\mathrm{s}\leq M-1                                                                                                                            \\\end{array}\right.
                                                                                                                                    \label{Eq_Lem_GRRWFsMNARsAWENOID_s_GRRWFs_001_001f}
\end{alignat}
Furthermore, for the values of $\{x,M_\pm,K_\mathrm{s}\}$ satisfying the conditions of \eqref{Eq_Lem_GRRWFsMNARsAWENOID_s_GRRWFs_001_001f},
\begin{alignat}{6}
\sum_{k_\mathrm{s}=0}^{K_\mathrm{s}}\sigma_{M_-,M_+,K_\mathrm{s},k_\mathrm{s}}(x)=1
                                                                                                                                    \label{Eq_Lem_GRRWFsMNARsAWENOID_s_GRRWFs_001_001g}
\end{alignat}
\end{subequations}
\end{lemma}
%
%
\begin{proof}
\begin{subequations}
                                                                                                                                    \label{Eq_Lem_GRRWFsMNARsAWENOID_s_GRRWFs_001_002}
Assume $M_\pm\in\mathbb{Z}\;:\;M:=M_-+M_+\geq3\Longrightarrow (M_--\ell_\mathrm{s})+(M_+-1+\ell_\mathrm{s})=M-1\geq2\;\forall\ell_\mathrm{s}\in\mathbb{Z}$.
Then, \eqref{Eq_Lem_GRRWFsMNARsAWENOID_s_GRRWFs_001_001d} applies to both functions $p_{M_--\ell_\mathrm{s},M_+-1+\ell_\mathrm{s}}$ ($\ell_\mathrm{ s}\in\{0,1\}$), and we have
\begin{alignat}{6}
p_{M_--\ell_\mathrm{s},M_+-1+\ell_\mathrm{s}}(x)\stackrel{\eqref{Eq_Lem_GRRWFsMNARsAWENOID_s_GRRWFs_001_001d}}{=}
\sum_{m_\mathrm{s}=0}^{1}\sigma_{M_--\ell_\mathrm{s},M_+-1+\ell_\mathrm{s},1,m_\mathrm{s}}(x)\;p_{(M_--\ell_\mathrm{s})-m_\mathrm{s},(M_+-1+\ell_\mathrm{ s})-1+m_\mathrm{s}}(x)
\qquad\left\{\begin{array}{l}\forall x\in{\mathbb R}\setminus\mathcal{S}_{\sigma_{M_--\ell_\mathrm{s},M_+-1+\ell_\mathrm{s},1}}\\
                             \forall M_\pm\in\mathbb{Z}\;:\;   M_-+M_+\geq3                                                     \\
                             \forall \ell_\mathrm{s}\in\{0,1\}                                                                  \\\end{array}\right.
                                                                                                                                    \label{Eq_Lem_GRRWFsMNARsAWENOID_s_GRRWFs_001_002a}
\end{alignat}
where $\sigma_{M_--\ell_\mathrm{s},M_+-1+\ell_\mathrm{s},1,m_\mathrm{s}}(x)$, being $1$-level weight-functions, exist by \eqrefsab{Eq_Lem_GRRWFsMNARsAWENOID_s_GRRWFs_001_001b}
                                                                                                                                  {Eq_Lem_GRRWFsMNARsAWENOID_s_GRRWFs_001_001d}.
Combining \eqrefsab{Eq_Lem_GRRWFsMNARsAWENOID_s_GRRWFs_001_001d}{Eq_Lem_GRRWFsMNARsAWENOID_s_GRRWFs_001_002a},
we have
\begin{alignat}{6}
p_{M_-,M_+}(x)\stackrel{\eqrefsab{Eq_Lem_GRRWFsMNARsAWENOID_s_GRRWFs_001_001d}
                                 {Eq_Lem_GRRWFsMNARsAWENOID_s_GRRWFs_001_002a}}{=}&
\sum_{\ell_\mathrm{s}=0}^{1}\sigma_{M_-,M_+,1,\ell_\mathrm{s}}(x)\left(\sum_{m_\mathrm{s}=0}^{1}\sigma_{M_--\ell_\mathrm{s},M_+-1+\ell_\mathrm{s},1,m_\mathrm{s}}(x)
                                                                                                \;p_{M_--\ell_\mathrm{s}-m_\mathrm{s},M_+-2+\ell_\mathrm{s}+m_\mathrm{s}}(x)\right)
                                                                                                                                    \notag\\
=&\sum_{\ell_\mathrm{s}=0}^{1}\sum_{m_\mathrm{s}=0}^{1}\sigma_{M_-,M_+,1,\ell_\mathrm{s}}(x)\;
                                                       \sigma_{M_--\ell_\mathrm{s},M_+-1+\ell_\mathrm{s},1,m_\mathrm{s}}(x)\;p_{M_--\ell_\mathrm{s}-m_\mathrm{s},M_+-2+\ell_\mathrm{s}+m_\mathrm{s}}(x)
                                                                                                                                    \notag\\
\stackrel{\text{\cite[(A.3)]{Gerolymos_2011a}}}{=}
 &\sum_{k_\mathrm{s}=0}^{2}\underbrace{\left(\sum_{\ell_\mathrm{s}=\max(0,k_\mathrm{s}-1)}^{\min(1,k_\mathrm{s})}
                                             \sigma_{M_-             ,M_+               ,1,          \ell_\mathrm{s}}(x)
                                           \;\sigma_{M_--\ell_\mathrm{s},M_+-1+\ell_\mathrm{s},1,k_\mathrm{s}-\ell_\mathrm{s}}(x)\right)}_{\displaystyle \sigma_{M_-,M_+,2,k_\mathrm{s}}(x)}
                                           \;p_{M_--k_\mathrm{s},M_+-2+k_\mathrm{s}}(x)
                                                                                                                                    \label{Eq_Lem_GRRWFsMNARsAWENOID_s_GRRWFs_001_002b}
\end{alignat}
\begin{alignat}{6}
\forall x\in{\mathbb R}\setminus\Big(\mathcal{S}_{\sigma_{M_-  ,M_+  ,1}}\cup
                                     \mathcal{S}_{\sigma_{M_-  ,M_+-1,1}}\cup
                                     \mathcal{S}_{\sigma_{M_--1,M_+  ,1}}\Big)\qquad\forall M_\pm\in\mathbb{Z}\;:\;M:=M_-+M_+\geq3\Longrightarrow 2\leq M-1
                                                                                                                                    \notag
\end{alignat}
which proves \eqrefsab{Eq_Lem_GRRWFsMNARsAWENOID_s_GRRWFs_001_001e}
                      {Eq_Lem_GRRWFsMNARsAWENOID_s_GRRWFs_001_001f},
for $K_\mathrm{s}=2$, because
\begin{alignat}{6}
                                \bigcup_{   L_\mathrm{s}=0}^{2           -1}\displaylimits
                                \bigcup_{\ell_\mathrm{s}=0}^{L_\mathrm{s}  }\mathcal{S}_{\sigma_{M_--\ell_\mathrm{s},M_+-L_\mathrm{s}+\ell_\mathrm{s},1}} 
                         = \Big(\bigcup_{\ell_\mathrm{s}=0}^{0             }\mathcal{S}_{\sigma_{M_--\ell_\mathrm{s},M_+             +\ell_\mathrm{s},1}}\Big)\cup
                           \Big(\bigcup_{\ell_\mathrm{s}=0}^{1             }\mathcal{S}_{\sigma_{M_--\ell_\mathrm{s},M_+-1           +\ell_\mathrm{s},1}}\Big)
                         = \Big(\mathcal{S}_{\sigma_{M_-  ,M_+  ,1}}\cup
                                \mathcal{S}_{\sigma_{M_-  ,M_+-1,1}}\cup
                                \mathcal{S}_{\sigma_{M_--1,M_+  ,1}}\Big)
                                                                                                                                    \label{Eq_Lem_GRRWFsMNARsAWENOID_s_GRRWFs_001_002c}
\end{alignat}
To prove \eqrefsab{Eq_Lem_GRRWFsMNARsAWENOID_s_GRRWFs_001_001e}
                  {Eq_Lem_GRRWFsMNARsAWENOID_s_GRRWFs_001_001f}
$\forall K_\mathrm{s}\in\{1,\cdots,M-1\}$, by induction, assume that \eqrefsab{Eq_Lem_GRRWFsMNARsAWENOID_s_GRRWFs_001_001e}
                                                                              {Eq_Lem_GRRWFsMNARsAWENOID_s_GRRWFs_001_001f}
are valid for some $K_\mathrm{s}-1\geq2$. Then
\begin{alignat}{6}
p_{M_-,M_+}(x)\stackrel{\eqref{Eq_Lem_GRRWFsMNARsAWENOID_s_GRRWFs_001_001f}}{=}\sum_{\ell_\mathrm{s}=0}^{K_\mathrm{s}-1}\sigma_{M_-,M_+,K_\mathrm{s}-1,\ell_\mathrm{s}}(x)\;p_{M_--\ell_\mathrm{s},M_+-(K_\mathrm{s}-1)+\ell_\mathrm{s}}(x)
\qquad\left\{\begin{array}{l}\displaystyle \forall x\in{\mathbb R}\setminus\bigcup_{   L_\mathrm{s}=0}^{K_\mathrm{s}-2}\displaylimits
                                                                           \bigcup_{\ell_\mathrm{s}=0}^{L_\mathrm{s}  }\mathcal{S}_{\sigma_{M_--\ell_\mathrm{s},M_+-L_\mathrm{s}+\ell_\mathrm{s},1}}\\
                                                                                                                                                                                                    \\
                                                                           \forall M_\pm\in\mathbb{Z}\;:\;M:=M_-+M_+\geq K_\mathrm{s}+1                                                             \\\end{array}\right.
                                                                                                                                    \label{Eq_Lem_GRRWFsMNARsAWENOID_s_GRRWFs_001_002d}
\end{alignat}
with $\sigma_{M_-,M_+,K_\mathrm{s}-1,\ell_\mathrm{s}}(x)$ in \eqref{Eq_Lem_GRRWFsMNARsAWENOID_s_GRRWFs_001_002d} defined by \eqref{Eq_Lem_GRRWFsMNARsAWENOID_s_GRRWFs_001_001e}.
Assuming $K_\mathrm{s}\leq M-1$ in \eqref{Eq_Lem_GRRWFsMNARsAWENOID_s_GRRWFs_001_002d},
we can further subdvide each of the $K_\mathrm{s}$ functions $p_{M_--\ell_\mathrm{s},M_+-(K_\mathrm{s}-1)+\ell_\mathrm{s}}(x)$ in \eqref{Eq_Lem_GRRWFsMNARsAWENOID_s_GRRWFs_001_002d},
once more, and we have by \eqref{Eq_Lem_GRRWFsMNARsAWENOID_s_GRRWFs_001_001d}
\begin{alignat}{6}
p_{M_--\ell_\mathrm{s},M_+-(K_\mathrm{s}-1)+\ell_\mathrm{s}}(x)\stackrel{\eqref{Eq_Lem_GRRWFsMNARsAWENOID_s_GRRWFs_001_001d}}{=}
\sum_{m_\mathrm{s}=0}^{1}\sigma_{M_--\ell_\mathrm{s},M_+-(K_\mathrm{s}-1)+\ell_\mathrm{s},1,m_\mathrm{s}}(x)\;p_{M_--\ell_\mathrm{s}-m_\mathrm{s},M_+-K_\mathrm{s}+\ell_\mathrm{s}+m_\mathrm{s}}(x)
\qquad\left\{\begin{array}{l}\displaystyle\forall x\in{\mathbb R}\setminus\mathcal{S}_{\sigma_{M_--\ell_\mathrm{s},M_+-(K_\mathrm{s}-1)+\ell_\mathrm{s},1}}\\
                                          \forall M_\pm\in\mathbb{Z}\;:\;M_-+M_+\geq K_\mathrm{s}+1                                                   \\
                                          \forall \ell_\mathrm{s}\in\{0,K_\mathrm{s}-1\}                                                                 \\\end{array}\right.
                                                                                                                                    \label{Eq_Lem_GRRWFsMNARsAWENOID_s_GRRWFs_001_002e}
\end{alignat}
where $\sigma_{M_--\ell_\mathrm{s},M_+-(K_\mathrm{s}-1)+\ell_\mathrm{s},1,m_\mathrm{s}}(x)$, being $1$-level weight-functions, exist by \eqrefsab{Eq_Lem_GRRWFsMNARsAWENOID_s_GRRWFs_001_001b}
                                                                                                                                                 {Eq_Lem_GRRWFsMNARsAWENOID_s_GRRWFs_001_001d}.
Combining \eqrefsab{Eq_Lem_GRRWFsMNARsAWENOID_s_GRRWFs_001_002d}
                   {Eq_Lem_GRRWFsMNARsAWENOID_s_GRRWFs_001_002e}, we have
\begin{alignat}{6}
p_{M_-,M_+}(x)\stackrel{\eqrefsab{Eq_Lem_GRRWFsMNARsAWENOID_s_GRRWFs_001_002d}
                                 {Eq_Lem_GRRWFsMNARsAWENOID_s_GRRWFs_001_002e}}{=}
 &\sum_{\ell_\mathrm{s}=0}^{K_\mathrm{s}-1}\sigma_{M_-,M_+,K_\mathrm{s}-1,\ell_\mathrm{s}}(x)
 \left(\sum_{m_\mathrm{s}=0}^{1}\sigma_{M_--\ell_\mathrm{s},M_+-(K_\mathrm{s}-1)+\ell_\mathrm{s},1,m_\mathrm{s}}(x)\;p_{M_--\ell_\mathrm{s}-m_\mathrm{s},M_+-K_\mathrm{s}+\ell_\mathrm{s}+m_\mathrm{s}}(x)\right)
                                                                                                                                    \notag\\
=&\sum_{\ell_\mathrm{s}=0}^{K_\mathrm{s}-1}\sum_{m_\mathrm{s}=0}^{1}
 \sigma_{M_-,M_+,K_\mathrm{s}-1,\ell_\mathrm{s}}(x)\;\sigma_{M_--\ell_\mathrm{s},M_+-(K_\mathrm{s}-1)+\ell_\mathrm{s},1,m_\mathrm{s}}(x)\;p_{M_--\ell_\mathrm{s}-m_\mathrm{s},M_+-K_\mathrm{s}+\ell_\mathrm{s}+m_\mathrm{s}}(x)
                                                                                                                                    \notag\\
\stackrel{\text{\cite[(A.3)]{Gerolymos_2011a}}}{=}
 &\sum_{k_\mathrm{s}=0}^{K_\mathrm{s}}\underbrace{\left(\sum_{\ell_\mathrm{s}=\max(0,k_\mathrm{s}-1)}^{\min(K_\mathrm{s}-1,k_\mathrm{s})}
                                                  \sigma_{M_-,M_+,K_\mathrm{s}-1,\ell_\mathrm{s}}(x)
                                                \;\sigma_{M_--\ell_\mathrm{s},M_+-(K_\mathrm{s}-1)+\ell_\mathrm{s},1,k_\mathrm{s}-\ell_\mathrm{s}}(x)\right)}_{\displaystyle \sigma_{M_-,M_+,K_\mathrm{s},k_\mathrm{s}}(x)}
                                                \;p_{M_--k_\mathrm{s},M_+-K_\mathrm{s}+k_\mathrm{s}}(x)
                                                                                                                                    \label{Eq_Lem_GRRWFsMNARsAWENOID_s_GRRWFs_001_002f}
\end{alignat}
\begin{alignat}{6}
\forall x\in{\mathbb R}\setminus\Bigg(\Big(\bigcup_{   L_\mathrm{s}=0}^{K_\mathrm{s}-2}\displaylimits
                                           \bigcup_{\ell_\mathrm{s}=0}^{L_\mathrm{s}  }\mathcal{S}_{\sigma_{M_--\ell_\mathrm{s},M_+-L_\mathrm{s}+\ell_\mathrm{s},1}}\Big)\cup
                                      \Big(\bigcup_{\ell_\mathrm{s}=0}^{K_\mathrm{s}-1}\displaylimits\mathcal{S}_{\sigma_{M_--\ell_\mathrm{s},M_+-(K_\mathrm{s}-1)+\ell_\mathrm{s},1}}\Big)\Bigg)
                                                                               \qquad \forall M_\pm\in\mathbb{Z}\;:\;M:=M_-+M_+\geq K_\mathrm{s}+1
                                                                                                                                    \notag
\end{alignat}
This completes the proof of \eqref{Eq_Lem_GRRWFsMNARsAWENOID_s_GRRWFs_001_001f} with weight-functions \eqref{Eq_Lem_GRRWFsMNARsAWENOID_s_GRRWFs_001_001e}, by induction.
By \eqref{Eq_Lem_GRRWFsMNARsAWENOID_s_GRRWFs_001_001e}, we have
\begin{alignat}{6}
 \sum_{k_\mathrm{s}=0}^{K_\mathrm{s}}\sigma_{M_-,M_+,K_\mathrm{s},k_\mathrm{s}}(x)
\stackrel{\eqref{Eq_Lem_GRRWFsMNARsAWENOID_s_GRRWFs_001_001e}}{=}&
 \sum_{k_\mathrm{s}=0}^{K_\mathrm{s}}\sum_{\ell_\mathrm{s}=\max(0,k_\mathrm{s}-1)}^{\min(K_\mathrm{s}-1,k_\mathrm{s})}\sigma_{M_-,M_+,K_\mathrm{s}-1,\ell_\mathrm{s}}(x)\;
                                                                                                                      \sigma_{M_--\ell_\mathrm{s},M_+-(K_\mathrm{s}-1)+\ell_\mathrm{s},1,k_\mathrm{s}-\ell_\mathrm{s}}(x)
                                                                                                                                    \notag\\
\stackrel{\text{\cite[(A.3)]{Gerolymos_2011a}}}{=}
 &\sum_{\ell_\mathrm{s}=0}^{K_\mathrm{s}-1}\sum_{m_\mathrm{s}=0}^{1}\sigma_{M_-,M_+,K_\mathrm{s}-1,\ell_\mathrm{s}}(x)\;
                                                                    \sigma_{M_--\ell_\mathrm{s},M_+-(K_\mathrm{s}-1)+\ell_\mathrm{s},1,m_\mathrm{s}}(x)
                                                                                                                                    \notag\\
=&
 \sum_{\ell_\mathrm{s}=0}^{K_\mathrm{s}-1}\left(\sigma_{M_-,M_+,K_\mathrm{s}-1,\ell_\mathrm{s}}(x)\;\sum_{m_\mathrm{s}=0}^{1}\sigma_{M_--\ell_\mathrm{s},M_+-(K_\mathrm{s}-1)+\ell_\mathrm{s},1,m_s}(x)\right)
\stackrel{\eqref{Eq_Lem_GRRWFsMNARsAWENOID_s_GRRWFs_001_001c}}{=}
 \sum_{\ell_\mathrm{s}=0}^{K_\mathrm{s}-1}\sigma_{M_-,M_+,K_\mathrm{s}-1,\ell_\mathrm{s}}(x)
                                                                                                                                    \label{Eq_Lem_GRRWFsMNARsAWENOID_s_GRRWFs_001_002g}
\end{alignat}
{\em ie} the sum of the weight-functions \eqref{Eq_Lem_GRRWFsMNARsAWENOID_s_GRRWFs_001_001e} is the same $\forall K_\mathrm{s}\in\{1,\cdots,M-1\}$ (by induction).
Since, by \eqref{Eq_Lem_GRRWFsMNARsAWENOID_s_GRRWFs_001_001c}, \eqref{Eq_Lem_GRRWFsMNARsAWENOID_s_GRRWFs_001_001g} holds for $K_\mathrm{s}=1$, \eqref{Eq_Lem_GRRWFsMNARsAWENOID_s_GRRWFs_001_002g}
suffices to prove \eqref{Eq_Lem_GRRWFsMNARsAWENOID_s_GRRWFs_001_001g} $\forall K_\mathrm{s}\in\{1,\cdots,M-1\}$, by induction.
\end{subequations}
\qed
\end{proof}
%

%
%
%
%
%
%
%
%
%
\section{Application to the Lagrange interpolating polynomial}\label{GRRWFsMNARsAWENOID_s_ALIP}
%
%
%
%
%
%
%
%
%

By Aitken's Lemma \cite[pp. 204--205]{Henrici_1964a}, the Lagrange interpolating polynomial
satisfies the $1$-level subdivision property \eqrefsab{Eq_Lem_GRRWFsMNARsAWENOID_s_GRRWFs_001_001b}
                                                      {Eq_Lem_GRRWFsMNARsAWENOID_s_GRRWFs_001_001c},
with weight-functions which are linear polynomials, and therefore defined $\forall x\in\mathbb{R}$,
implying that $\mathcal{S}_{I,\tsc{x}_{i-M_-,i+M_+},K_\mathrm{s}}=\varnothing$ in \eqref{Eq_Lem_GRRWFsMNARsAWENOID_s_GRRWFs_001_001}.
Application of \lemrefnp{Lem_GRRWFsMNARsAWENOID_s_GRRWFs_001} to the Lagrange interpolating polynomial can be summarized in the following proposition.
%
\begin{proposition}[{\rm Weight-functions for the Lagrange interpolating polynomial}]
\label{Prp_GRRWFsMNARsAWENOID_s_ALIP_001}
\begin{subequations}
                                                                                                       \label{Eq_Prp_GRRWFsMNARsAWENOID_s_ALIP_001_001}
Assume the conditions of \defrefnp{Def_GRRWFsMNARsAWENOID_s_I_002}. Then, the weight-functions $\sigma_{I,\tsc{x}_{i-M_-,i+M_+},K_\mathrm{s},k_\mathrm{s}}(x)$
in the representation \eqref{Eq_GRRWFsMNARsAWENOID_s_I_001} of the Lagrange interpolating polynomial on $\tsc{x}_{i-M_-,i+M_+}$, $p_{I,\tsc{x}_{i-M_-,i+M_+}}(x;f)$,
are real polynomials of degree $K_\mathrm{s}$ with only real roots, expressed by
\begin{alignat}{6}
&\mathbb{R}_{K_\mathrm{s}}[x]\ni\sigma_{I,\tsc{x}_{i-M_-,i+M_+},K_\mathrm{s},k_\mathrm{s}}(x):=(-1)^{K_\mathrm{s}-k_\mathrm{s}}\ss_{I,\tsc{x}_{i-M_-,i+M_+},K_\mathrm{s},k_\mathrm{s}}
                                                                                               \prod_{x_n\in\tsc{x}_{i-M_-,i+M_+}\setminus\tsc{x}_{i-M_-+k_\mathrm{s},i+M_+-K_\mathrm{s}+k_\mathrm{s}}}(x-x_n)&
                                                                                                                                    \label{Eq_Prp_GRRWFsMNARsAWENOID_s_ALIP_001_001a}\\
&\forall x\in\mathbb{R}\qquad\forall k_\mathrm{s}\in\{0,\cdots,K_\mathrm{s}\}\qquad\forall K_\mathrm{s}\in\{1,\cdots,M-1:=M_-+M_+-1\}&
                                                                                                                                    \notag
\end{alignat}
where the strictly positive real numbers $\ss_{I,\tsc{x}_{i-M_-,i+M_+},K_\mathrm{s},k_\mathrm{s}}$ depend on the points of the stencil $\tsc{x}_{i-M_-,i+M_+}$ \defref{Def_GRRWFsMNARsAWENOID_s_I_001},
and are generated by the recurrence relation
\begin{alignat}{6}
&\mathbb{R}_{>0}\ni \ss_{I,\tsc{x}_{i-M_-,i+M_+},K_\mathrm{s},k_\mathrm{s}}:=\left\{\begin{array}{ll}\dfrac{1}{x_{i+M_+}-x_{i-M_-}}                                  &K_\mathrm{s}=1       \\
                                                                                                                                                                   &                     \\
                                                                                                   {\displaystyle\sum_{\ell_\mathrm{s}=\max(0,k_\mathrm{s}-1)}^{\min(K_\mathrm{s}-1,k_\mathrm{s})}
                                                                                                   \ss_{I,\tsc{x}_{i-M_-                ,i+M_+                                 },K_\mathrm{s}-1,             \ell_\mathrm{s}}\;
                                                                                                   \ss_{I,\tsc{x}_{i-M_-+\ell_\mathrm{s},i+M_+-(K_\mathrm{s}-1)+\ell_\mathrm{s}},             1,k_\mathrm{s}-\ell_\mathrm{s}}
                                                                             }                                                                                     &K_\mathrm{s}\geq2    \\\end{array}\right.
                                                                                                                                    \notag\\
&\forall k_\mathrm{s}\in\{0,\cdots,K_\mathrm{s}\}\qquad\forall K_\mathrm{s}\in\{1,\cdots,M-1:=M_-+M_+-1\}
                                                                                                       \label{Eq_Prp_GRRWFsMNARsAWENOID_s_ALIP_001_001b}
\end{alignat}
The weight-functions \eqref{Eq_Prp_GRRWFsMNARsAWENOID_s_ALIP_001_001a} satisfy the consistency condition \eqref{Eq_GRRWFsMNARsAWENOID_s_I_001b} and the recurrence relation \eqref{Eq_Lem_GRRWFsMNARsAWENOID_s_GRRWFs_001_001e}.
\end{subequations}
\end{proposition}
%
%
\begin{proof}
\begin{subequations}
                                                                                                       \label{Eq_Prp_GRRWFsMNARsAWENOID_s_ALIP_001_002}
The case $K_{\rm s}=1$
\begin{alignat}{6}
\sigma_{I,\tsc{x}_{i-M_-,i+M_+},1,k_\mathrm{s}}(x)\stackrel{\eqrefsab{Eq_Prp_GRRWFsMNARsAWENOID_s_ALIP_001_001a}
                                                                     {Eq_Prp_GRRWFsMNARsAWENOID_s_ALIP_001_001b}}{=}(-1)^{1-k_\mathrm{s}}\dfrac{1}{x_{i+M_+}-x_{i-M_-}}(x-x_{i+M_+-k_\mathrm{s}M})
\quad\forall k_\mathrm{s}\in\{0,1\}
                                                                                                       \label{Eq_Prp_GRRWFsMNARsAWENOID_s_ALIP_001_002a}
\end{alignat}
holds since it is exactly Aitken's Lemma \cite[pp. 204--205]{Henrici_1964a}. Since \eqref{Eq_Prp_GRRWFsMNARsAWENOID_s_ALIP_001_001} hold for $K_\mathrm{s}=1$, $\forall M_\pm\in\mathbb{Z}: M:=M_-+M_+\geq2$
the family of Lagrange interpolating polynomials is equipped with the $1$-level subdivison rule \eqrefsatob{Eq_Lem_GRRWFsMNARsAWENOID_s_GRRWFs_001_001a}
                                                                                                           {Eq_Lem_GRRWFsMNARsAWENOID_s_GRRWFs_001_001c},
and therefore satisfies the conditions of \lemrefnp{Lem_GRRWFsMNARsAWENOID_s_GRRWFs_001}.
We can therefore apply \eqref{Eq_Lem_GRRWFsMNARsAWENOID_s_GRRWFs_001_001e} to $\sigma_{I,\tsc{x}_{i-M_-,i+M_+},K_\mathrm{s},k_\mathrm{s}}(x)$.
To obtain the simpler expressions \eqref{Eq_Prp_GRRWFsMNARsAWENOID_s_ALIP_001_001}, assume that \eqref{Eq_Prp_GRRWFsMNARsAWENOID_s_ALIP_001_001a}
holds for some $K_\mathrm{s}-1\geq1$. Then by \lemrefnp{Lem_GRRWFsMNARsAWENOID_s_GRRWFs_001}
\begin{alignat}{6}
\sigma_{I,\tsc{x}_{i-M_-,i+M_+},K_\mathrm{s},k_\mathrm{s}}(x)\stackrel{\eqref{Eq_Lem_GRRWFsMNARsAWENOID_s_GRRWFs_001_001e}}{=}&
\sum_{\ell_\mathrm{s}=\max(0,k_\mathrm{s}-1)}^{\min(K_\mathrm{s}-1,k_\mathrm{s})}
\sigma_{I,\tsc{x}_{i-M_-                ,i+M_+                                 },K_\mathrm{s}-1,             \ell_\mathrm{s}}(x)\;
\sigma_{I,\tsc{x}_{i-M_-+\ell_\mathrm{s},i+M_+-(K_\mathrm{s}-1)+\ell_\mathrm{s}},             1,k_\mathrm{s}-\ell_\mathrm{s}}(x)
                                                                                                                                    \notag\\
\stackrel{\eqref{Eq_Prp_GRRWFsMNARsAWENOID_s_ALIP_001_001}}{=}&\sum_{\ell_\mathrm{s}=\max(0,k_\mathrm{s}-1)}^{\min(K_\mathrm{s}-1,k_\mathrm{s})}\Bigg(
                                                              (-1)^{K_\mathrm{s}-1-\ell_\mathrm{s}}\ss_{I,\tsc{x}_{i-M_-,i+M_+},K_\mathrm{s}-1,\ell_\mathrm{s}}
                                                              \prod_{x_n\in\tsc{x}_{i-M_-,i+M_+}\setminus\tsc{x}_{i-M_-+\ell_\mathrm{s},i+M_+-(K_\mathrm{s}-1)+\ell_\mathrm{s}}}(x-x_n)
                                                                                                                                    \notag\\
                                                             &\times(-1)^{1-k_\mathrm{s}+\ell_\mathrm{s}}\ss_{I,\tsc{x}_{i-M_-+\ell_\mathrm{s},i+M_+-(K_\mathrm{s}-1)+\ell_\mathrm{s}},1,k_\mathrm{s}-\ell_\mathrm{s}}
                                                              \prod_{x_n\in\tsc{x}_{i-M_-+\ell_\mathrm{s},i+M_+-(K_\mathrm{s}-1)+\ell_\mathrm{s}}\setminus
                                                                           \tsc{x}_{i-M_-+\ell_\mathrm{s}+k_\mathrm{s}-\ell_\mathrm{s},i+M_+-(K_\mathrm{s}-1)+\ell_\mathrm{s}-1+k_\mathrm{s}-\ell_\mathrm{s}}}(x-x_n)
                                                                                                                                              \Bigg)
                                                                                                                                    \notag\\
                                                            =&(-1)^{K_\mathrm{s}-k_\mathrm{s}}\left(\sum_{\ell_\mathrm{s}=\max(0,k_\mathrm{s}-1)}^{\min(K_\mathrm{s}-1,k_\mathrm{s})}
                                                                                               \ss_{I,\tsc{x}_{i-M_-,i+M_+},K_\mathrm{s}-1,\ell_\mathrm{s}}
                                                                                               \ss_{I,\tsc{x}_{i-M_-+\ell_\mathrm{s},i+M_+-(K_\mathrm{s}-1)+\ell_\mathrm{s}},1,k_\mathrm{s}-\ell_\mathrm{s}}\right)
                                                                                               \prod_{x_n\in\tsc{x}_{i-M_-,i+M_+}\setminus\tsc{x}_{i-M_-+k_\mathrm{s},i+M_+-K_\mathrm{s}+k_\mathrm{s}}}(x-x_n)
                                                                                                       \label{Eq_Prp_GRRWFsMNARsAWENOID_s_ALIP_001_002b}
\end{alignat}
because $\Big(\tsc{x}_{i-M_-,i+M_+}\setminus\tsc{x}_{i-M_-+\ell_\mathrm{s},i+M_+-(K_\mathrm{s}-1)+\ell_\mathrm{s}}\Big)\cup
         \Big(\tsc{x}_{i-M_-+\ell_\mathrm{s},i+M_+-(K_\mathrm{s}-1)+\ell_\mathrm{s}}\setminus\tsc{x}_{i-M_-+\ell_\mathrm{s}+k_\mathrm{s}-\ell_\mathrm{s},i+M_+-(K_\mathrm{s}-1)-1+\ell_\mathrm{s}+k_\mathrm{s}-\ell_\mathrm{s}}\Big)=
         \Big(\tsc{x}_{i-M_-,i+M_+}\setminus\tsc{x}_{i-M_-+k_\mathrm{s},i+M_+-K_\mathrm{s}+k_\mathrm{s}}\Big)$.
Since \eqrefsab{Eq_Prp_GRRWFsMNARsAWENOID_s_ALIP_001_001a}
               {Eq_Prp_GRRWFsMNARsAWENOID_s_ALIP_001_001b} are valid for $K_\mathrm{s}=1$ by Aitken's Lemma \cite[pp. 204--205]{Henrici_1964a},
\eqref{Eq_Prp_GRRWFsMNARsAWENOID_s_ALIP_001_002b} proves that they are also valid for $K_\mathrm{s}=2$, and by induction $\forall K_\mathrm{s}\in\{1,\cdots,M-1:=M_-+M_+-1\}$,
completing the proof. Direct computation, using \eqref{Eq_Prp_GRRWFsMNARsAWENOID_s_ALIP_001_002a}, proves that the consistency relation \eqref{Eq_GRRWFsMNARsAWENOID_s_I_001b} holds for $K_\mathrm{s}=1$,
and, by \lemrefnp{Lem_GRRWFsMNARsAWENOID_s_GRRWFs_001}, $\forall K_\mathrm{s}\in\{1,\cdots,M-1\}$. Finally, strict positivity of
$\ss_{I,\tsc{x}_{i-M_-,i+M_+},1,0}\stackrel{\eqref{Eq_Prp_GRRWFsMNARsAWENOID_s_ALIP_001_001b}}{=}\ss_{I,\tsc{x}_{i-M_-,i+M_+},1,1}$ follows by the order relations assumed in \eqref{Eq_Def_GRRWFsMNARsAWENOID_s_I_002_001a},
and then by induction, using \eqref{Eq_Prp_GRRWFsMNARsAWENOID_s_ALIP_001_001b}, $\forall k_\mathrm{s}\in\{0,\cdots,K_\mathrm{s}\}$ and $\forall K_\mathrm{s}\in\{1,\cdots,M-1:=M_-+M_+-1\}$.
\end{subequations}
\qed
\end{proof}
%

Because of the positivity of the numbers $\ss_{I,\tsc{x}_{i-M_-,i+M_+},K_\mathrm{s},k_\mathrm{s}}\in\mathbb{R}_{>0}$ \eqref{Eq_Prp_GRRWFsMNARsAWENOID_s_ALIP_001_001b}
it is quite straightforward to study the sign of the weight-functions $\sigma_{I,\tsc{x}_{i-M_-,i+M_+},K_\mathrm{s},k_\mathrm{s}}(x)$ \eqref{Eq_Prp_GRRWFsMNARsAWENOID_s_ALIP_001_001a},
which allows to determine the intervals on the real axis where the combination \eqref{Eq_GRRWFsMNARsAWENOID_s_I_001} of the Lagrange interpolating polynomials on the substenicls is convex.
%
\begin{proposition}[{\rm Convexity in the neighbourhood of $x_i$}]
\label{Prp_GRRWFsMNARsAWENOID_s_ALIP_002}
Assume the conditions of \prprefnp{Prp_GRRWFsMNARsAWENOID_s_ALIP_001}. Furthermore assume that $K_\mathrm{s}\leq\left\lceil\frac{M}{2}\right\rceil$.
Then the weight-functions of the combination \eqref{Eq_GRRWFsMNARsAWENOID_s_I_001} of the Lagrange interpolating polynomials on substencils \prpref{Prp_GRRWFsMNARsAWENOID_s_ALIP_001}
satisfy
\begin{alignat}{6}
0\leq\sigma_{I,\tsc{x}_{i-M_-,i+M_+},K_\mathrm{s},k_\mathrm{s}}(x)\leq1\qquad\forall x\in[x_{i-M_-+K_\mathrm{s}-1},
                                                                          x_{i+M_+-K_\mathrm{s}+1}]\qquad\left\{\begin{array}{l}\forall M_\pm\in\mathbb{Z}:M:=M_-+M_+\geq2                            \\
                                                                                                                                \forall K_\mathrm{s}\in\{1,\cdots,\left\lceil\frac{M}{2}\right\rceil\}\\
                                                                                                                                \forall k_\mathrm{s}\in\{0,\cdots,K_\mathrm{s}\}                      \\\end{array}\right.
                                                                                                       \label{Eq_Prp_GRRWFsMNARsAWENOID_s_ALIP_002_001}
\end{alignat}
\end{proposition}
%
%
\begin{proof}
\begin{subequations}
                                                                                                       \label{Eq_Prp_GRRWFsMNARsAWENOID_s_ALIP_002_002}
Because of the consistency condition \eqref{Eq_GRRWFsMNARsAWENOID_s_I_001b},
(non strict) positivity of the weight-functions $\sigma_{I,\tsc{x}_{i-M_-,i+M_+},K_\mathrm{s},k_\mathrm{s}}(x)$ \eqref{Eq_Prp_GRRWFsMNARsAWENOID_s_ALIP_001_001}
suffices (proof by contradiction) to prove \eqref{Eq_Prp_GRRWFsMNARsAWENOID_s_ALIP_002_001}.
Rewrite \eqref{Eq_Prp_GRRWFsMNARsAWENOID_s_ALIP_001_001a} as
\begin{alignat}{6}
\sigma_{I,\tsc{x}_{i-M_-,i+M_+},K_\mathrm{s},k_\mathrm{s}}(x)\stackrel{\eqref{Eq_Prp_GRRWFsMNARsAWENOID_s_ALIP_001_001a}}{=}
(-1)^{K_\mathrm{s}-k_\mathrm{s}}\;\ss_{I,\tsc{x}_{i-M_-,i+M_+},K_\mathrm{s},k_\mathrm{s}}
\left\{\begin{array}{llcr}                                                          &\displaystyle\prod_{n=i+M_+-K_\mathrm{s}+k_\mathrm{s}+1}^{M_+}(x-x_n)    &     ;\;&   k_\mathrm{s}=0\\
                                                                                    &                                                                         &        &                            \\
                          \displaystyle\prod_{n=i-M_-}^{i-M_-+k_\mathrm{s}-1}(x-x_n)&\displaystyle\prod_{n=i+M_+-K_\mathrm{s}+k_\mathrm{s}+1}^{i+M_+}(x-x_n)  &     ;\;& 0<k_\mathrm{s}<K_\mathrm{s}\\
                                                                                    &                                                                         &        &                            \\
                          \displaystyle\prod_{n=i-M_-}^{i-M_-+k_\mathrm{s}-1}(x-x_n)&                                                                         &     ;\;&   k_\mathrm{s}=K_\mathrm{s} \\\end{array}\right.
                                                                                                       \label{Eq_Prp_GRRWFsMNARsAWENOID_s_ALIP_002_002a}
\end{alignat}
Obviously we have
\begin{alignat}{6}
\mathrm{sign}\left(\prod_{n=i-M_-}^{i-M_-+k_\mathrm{s}-1}(x-x_n)\right)=&1
\;\forall k_\mathrm{s}\in\{1,\cdots,K_\mathrm{s}\}\quad\forall x>\max_{0<k_\mathrm{s}\leq K_\mathrm{s}}x_{i-M_-+k_\mathrm{s}-1}\stackrel{\eqref{Eq_Def_GRRWFsMNARsAWENOID_s_I_001_001}}{=}x_{i-M_-+K_\mathrm{s}-1}
                                                                                                       \label{Eq_Prp_GRRWFsMNARsAWENOID_s_ALIP_002_002b}\\
\mathrm{sign}\left(\prod_{n=i+M_+-K_\mathrm{s}+k_\mathrm{s}+1}^{i+M_+}(x-x_n)\right)=&(-1)^{K_\mathrm{s}-k_\mathrm{s}}
\;\forall k_\mathrm{s}\in\{0,\cdots,K_\mathrm{s}-1\}\quad\forall x<\min_{0\leq k_\mathrm{s}< K_\mathrm{s}}x_{i+M_+-K_\mathrm{s}+k_\mathrm{s}+1}\stackrel{\eqref{Eq_Def_GRRWFsMNARsAWENOID_s_I_001_001}}{=}x_{i+M_+-K_\mathrm{s}+1}
                                                                                                       \label{Eq_Prp_GRRWFsMNARsAWENOID_s_ALIP_002_002c}
\end{alignat}
Combining \eqrefsatob{Eq_Prp_GRRWFsMNARsAWENOID_s_ALIP_002_002a}{Eq_Prp_GRRWFsMNARsAWENOID_s_ALIP_002_002c} with the positivity of the numbers
$\ss_{I,\tsc{x}_{i-M_-,i+M_+},K_\mathrm{s},k_\mathrm{s}}\in\mathbb{R}_{>0}$ \eqref{Eq_Prp_GRRWFsMNARsAWENOID_s_ALIP_001_001b}, and taking into account that $(-1)^{K_\mathrm{s}-K_\mathrm{s}}=1$,
proves \eqref{Eq_Prp_GRRWFsMNARsAWENOID_s_ALIP_002_001}. Notice that the condition for the interval $[x_{i-M_-+K_\mathrm{s}-1},x_{i+M_+-K_\mathrm{s}+1}]$ in \eqref{Eq_Prp_GRRWFsMNARsAWENOID_s_ALIP_002_001}
to contain at least 1 cell (at least 2 grid-points) is
$-M_-+K_\mathrm{s}-1< M_+-K_\mathrm{s}+1 \iff 2K_\mathrm{s} < M_++M_-+2 \stackrel{\text{\cite[(A.2)]{Gerolymos_2011a}}}{\iff} K_\mathrm{s} < \left\lceil\frac{M+2}{2}\right\rceil=\left\lceil\frac{M}{2}\right\rceil+1$,
which explains the additional constraint on $K_\mathrm{s}$ included in the hypotheses of \prprefnp{Prp_GRRWFsMNARsAWENOID_s_ALIP_002}.
\end{subequations}
\qed
\end{proof}
%
%
\begin{remark}[{\rm Typical stencils \cite{Carlini_Ferretti_Russo_2005a,Liu_Shu_Zhang_2009a}}]
\label{Rmk_GRRWFsMNARsAWENOID_s_GRRWFs_001}
For $\sigma_{I,\tsc{x}_{i-(r-1),i+r},r-1,k_\mathrm{s}}(x)$ the positivity interval is, by \eqref{Eq_Prp_GRRWFsMNARsAWENOID_s_ALIP_002_001},
$[x_{i-(r-1)+(r-1)-1}, x_{i+r-(r-1)+1}]=[x_{i-1},x_{i+2}]$ in agreement with \cite[Tab. 2.1, p. 507]{Liu_Shu_Zhang_2009a},
while for $\sigma_{I,\tsc{x}_{i-r,i+r},r,k_\mathrm{s}}(x)$ the positivity interval is, by \eqref{Eq_Prp_GRRWFsMNARsAWENOID_s_ALIP_002_001},
$[x_{i-r+r-1}, x_{i+r-r+1}]=[x_{i-1},x_{i+1}]$ in agreement with \cite[Tab. 2.2, p. 511]{Liu_Shu_Zhang_2009a}.
\qed
\end{remark}
%

\lemrefnp{Lem_GRRWFsMNARsAWENOID_s_GRRWFs_001} only requires the determination of weight-functions for the $1$-level subdivision.
It is therefore not limited to a particular family of stencils and/or subdivisions, and can be used to determine weight-functions
on biased stencils, {\em eg} near the boundaries of the computational domain.

%
\begin{remark}[{\rm Relation to previous work}]
\label{Rmk_GRRWFsMNARsAWENOID_s_GRRWFs_002}
\tsc{weno} interpolation applied to the development of central \tsc{weno} schemes only requires the computation of the value of the weight-functions
at specific points on the stencil, and these can be computed by solving a linear system \cite{Shu_2009a}.
Carlini \etal~\cite{Carlini_Ferretti_Russo_2005a} pointed out that the weight-functions in representation \eqref{Eq_GRRWFsMNARsAWENOID_s_I_001}
are of the form \eqref{Eq_Prp_GRRWFsMNARsAWENOID_s_ALIP_001_001a}, with unknown constants $\gamma_{k_\mathrm{s}}:=(-1)^{K_\mathrm{s}-k_\mathrm{s}}\ss_{I,\tsc{x}_{i-M_-,i+M_+},K_\mathrm{s},k_\mathrm{s}}$,
which can be determined, in the general case, by solution of a linear triangular system. These authors \cite{Carlini_Ferretti_Russo_2005a} studied in particular the stencil
$\tsc{x}_{i-(r-1),i+r}$, for $K_\mathrm{s}=r-1$, for which they obtained an analytical expression for the coefficients $\gamma_{k_\mathrm{s}}$,
and proved convexity $\forall x\in[x_i,x_{i+1}]$. Liu \etal~\cite{Liu_Shu_Zhang_2009a} used the same form as Carlini \etal~\cite{Carlini_Ferretti_Russo_2005a} for the weight-functions,
computed the coefficients up to $r=6$, and observed that the interval of convexity is actually $[x_{i-1},x_{i+2}]$ \rmkref{Rmk_GRRWFsMNARsAWENOID_s_GRRWFs_001}.
They also studied the stencil $\tsc{x}_{i-r,i+r}$, for $K_\mathrm{s}=r$, gave an analytical expression for the coefficients $\gamma_{k_\mathrm{s}}$, which were computed
up to $r=5$, and observed that the interval of convexity in this case is $[x_{i-1},x_{i+1}]$ \rmkref{Rmk_GRRWFsMNARsAWENOID_s_GRRWFs_001}.
These results are in agreement with those proven in \prprefnp{Prp_GRRWFsMNARsAWENOID_s_ALIP_002} \rmkref{Rmk_GRRWFsMNARsAWENOID_s_GRRWFs_001}.
\prprefnp{Prp_GRRWFsMNARsAWENOID_s_ALIP_001} studies an arbitrary stencil $\tsc{x}_{i-M_-,i+M_+}$ \defref{Def_GRRWFsMNARsAWENOID_s_I_001}, and level of subdivision $K_\mathrm{s}\in\{1,\cdots, M-1\}$ ($M:=M_-+M_+$),
and obtains an analytical recursive expression for the coefficients $\ss_{I,\tsc{x}_{i-M_-,i+M_+},K_\mathrm{s},k_\mathrm{s}}$ in \eqref{Eq_Prp_GRRWFsMNARsAWENOID_s_ALIP_001_001a}.
In this way, we were able to give a formal proof for the interval of convexity, which was determined for general values of $M_\pm$ and $K_\mathrm{s}$ \prpref{Prp_GRRWFsMNARsAWENOID_s_ALIP_002}.
\qed
\end{remark}
%

%
\begin{remark}[{\rm Alternative formulation}]
\label{Rmk_GRRWFsMNARsAWENOID_s_GRRWFs_003}
The expression of the $(K_\mathrm{s}=1)$-level weight-functions \eqref{Eq_Prp_GRRWFsMNARsAWENOID_s_ALIP_001_002a} can also be written in an equivalent form,
using ratios of fundamental functions of Lagrange interpolation on $\tsc{x}_{i-M_-,i+M_+}$ \defref{Def_GRRWFsMNARsAWENOID_s_I_001} and on its
$(K_\mathrm{s}=1)$-level Neville substencils \defref{Def_GRRWFsMNARsAWENOID_s_I_002}, $\tsc{x}_{i-M_-,i+M_+-1}$ and $\tsc{x}_{i-M_-+1,i+M_+}$.
The Lagrange interpolating polynomial of a function $f:\mathbb{R}\to\mathbb{R}$, on $\tsc{x}_{i-M_-,i+M_+}$, can be expressed \cite[(9.3,9.4), p. 184]{Henrici_1964a} as
\begin{subequations}
                                                                                                       \label{Eq_Rmk_GRRWFsMNARsAWENOID_s_GRRWFs_003_001}
\begin{alignat}{6}
p_{I,\tsc{x}_{i-M_-,i+M_+}}(x;f)\stackrel{\text{\cite[(9.3)]{Henrici_1964a}}}{=}&\sum_{\ell=-M_-}^{M_+}\alpha_{I,\tsc{x}_{i-M_-,i+M_+},i+\ell}(x)\;f_{i+\ell}
                                                                                                       \label{Eq_Rmk_GRRWFsMNARsAWENOID_s_GRRWFs_003_001a}\\
\alpha_{I,\tsc{x}_{i-M_-,i+M_+},i+\ell}(x)\stackrel{\text{\cite[(9.4)]{Henrici_1964a}}}{:=}&\prod_{\substack{k=-M_-\\
                                                                                                            k\neq\ell}}^{M_+}{\dfrac{x         -x_{i+k}}
                                                                                                                                    {x_{i+\ell}-x_{i+k}}}
&\forall\;\ell\in\{-M_-,\cdots,M_+\}
                                                                                                       \label{Eq_Rmk_GRRWFsMNARsAWENOID_s_GRRWFs_003_001b}\\
f_{i+\ell}:=&f(x_{i+\ell})
&\forall\;\ell\in\{-M_-,\cdots,M_+\}
                                                                                                       \label{Eq_Rmk_GRRWFsMNARsAWENOID_s_GRRWFs_003_001c}
\end{alignat}
\end{subequations}
where the $M+1$ polynomials $\alpha_{I,\tsc{x}_{i-M_-,i+M_+},i+\ell}\in\mathbb{R}_M[x]$ \eqref{Eq_Rmk_GRRWFsMNARsAWENOID_s_GRRWFs_003_001b} are $\neq0_{\mathbb{R}_M[x]}(x)$,
linearly independent, and form a basis of the space of all polynomials with real coefficients and degree $\leq M$, $\mathbb{R}_M[x]$ \cite[p. 2771]{Gerolymos_2012a}.
It can be verified by direct computation, using definition \eqref{Eq_Rmk_GRRWFsMNARsAWENOID_s_GRRWFs_003_001b}, that \eqref{Eq_Prp_GRRWFsMNARsAWENOID_s_ALIP_001_002a} is equivalent to
\begin{subequations}
                                                                                                       \label{Eq_Rmk_GRRWFsMNARsAWENOID_s_GRRWFs_003_002}
\begin{alignat}{6}
\sigma_{I,\tsc{x}_{i-M_-,i+M_+},1,0}(x)\stackrel{\eqrefsab{Eq_Prp_GRRWFsMNARsAWENOID_s_ALIP_001_002a}
                                                          {Eq_Rmk_GRRWFsMNARsAWENOID_s_GRRWFs_003_001b}}{=}\dfrac{\alpha_{I,\tsc{x}_{i-M_-,i+M_+},i-M_-}(x)}
                                                                                                                 {\alpha_{I,\tsc{x}_{i-M_-,i+M_+-1},i-M_-}(x)}
                                                                                                       \label{Eq_Rmk_GRRWFsMNARsAWENOID_s_GRRWFs_003_002a}\\
\sigma_{I,\tsc{x}_{i-M_-,i+M_+},1,1}(x)\stackrel{\eqrefsab{Eq_Prp_GRRWFsMNARsAWENOID_s_ALIP_001_002a}
                                                          {Eq_Rmk_GRRWFsMNARsAWENOID_s_GRRWFs_003_001b}}{=}\dfrac{\alpha_{I,\tsc{x}_{i-M_-,i+M_+},i+M_+}(x)}
                                                                                                                 {\alpha_{I,\tsc{x}_{i-M_-+1,i+M_+},i+M_+}(x)}
                                                                                                       \label{Eq_Rmk_GRRWFsMNARsAWENOID_s_GRRWFs_003_002b}
\end{alignat}
\end{subequations}
The interest of this alternative expression \eqref{Eq_Rmk_GRRWFsMNARsAWENOID_s_GRRWFs_003_002}, which is analogous to the expression of the $(K_\mathrm{s}=1)$-level weight-functions for the representation of the
Lagrange reconstructing polynomial on a homogeneous grid \cite[Lemma 4.2, p. 2780]{Gerolymos_2012a}, is that it can be generalized for the representation
of the $n$-derivative of the Lagrange interpolating polynomial by combination of substencils, as will be shown in \parrefnp{GRRWFsMNARsAWENOID_s_AnDLIP}.
\qed
\end{remark}
%

%
%
%
%
%
%
%
%
%
\section{Application to the $n$-derivative of the Lagrange intepolating polynomial}\label{GRRWFsMNARsAWENOID_s_AnDLIP}
%
%
%
%
%
%
%
%
%

One of the motivations that led to the formulation of \lemrefnp{Lem_GRRWFsMNARsAWENOID_s_GRRWFs_001} was the study of \tsc{weno} reconstruction in view of the computation of $f'(x)$,
and this application is studied in \cite{Gerolymos_2012a} (results and relation to previous work are summarized in \rmkrefnp{Rmk_GRRWFsMNARsAWENOID_s_AnDLIP_001}).
The expression of the $(K_\mathrm{s}=1)$-level weight-functions for the Lagrange reconstructing polynomial is similar to \eqref{Eq_Rmk_GRRWFsMNARsAWENOID_s_GRRWFs_003_002},
upon replacing the fundamental functions of Lagrange interpolation in \eqref{Eq_Rmk_GRRWFsMNARsAWENOID_s_GRRWFs_003_002} by the
corresponding fundamental functions of Lagrange reconstruction \cite[(32), Lemma 4.2, pp. 2780--2781]{Gerolymos_2012a}. It turns out that a similar
relation is valid for the $n$-derivative of the Lagrange interpolating polynomial
\begin{alignat}{6}
p^{(n)}_{I,\tsc{x}_{i-M_-,i+M_+}}(x;f)\stackrel{\eqref{Eq_Rmk_GRRWFsMNARsAWENOID_s_GRRWFs_003_001a}}{=}&\sum_{\ell=-M_-}^{M_+}\alpha^{(n)}_{I,\tsc{x}_{i-M_-,i+M_+},i+\ell}(x)\;f_{i+\ell}
                                                                                                       \label{Eq_GRRWFsMNARsAWENOID_s_AnDLIP_001}
\end{alignat}
and these new results are formulated in \prprefnp{Prp_GRRWFsMNARsAWENOID_s_AnDLIP_001}.

%
\begin{remark}[{\rm Lagrange reconstructing polynomial}]
\label{Rmk_GRRWFsMNARsAWENOID_s_AnDLIP_001}
The case of reconstruction of a function $h(x)$ from its cell-averages $f(x)$, sampled on a given stencil \defref{Def_GRRWFsMNARsAWENOID_s_I_001}, is important for the construction of numerical schemes used in the
solution of hyperbolic \tsc{pde}s \cite{Shu_2009a}. In the particular case of homogeneous grids ($x_{i+1}-x_{i}=\Delta x=\const\in\mathbb{R}_{>0}\;\forall i$),
reconstruction can be used for the computation of $f'(x)$ \cite{Shu_2009a,
                                                                Gerolymos_2011a}.\footnote{\label{ff_GRRWFsMNARsAWENOID_s_EoA_ss_R1_001}
                                                                                           on general inhomogeneous grids reconstruction does not provide $f'(x)$,
                                                                                           only numerical fluxes for the discretization of the \tsc{pde} \cite{Shu_2009a}
                                                                                          }
In the context of methods for the determination of numerical fluxes, what is needed are the values of the weight-functions at $x_{i+\frac{1}{2}}$ (optimal weights [3]),
which were usually computed from the solution of a linear system \cite[(13), p. 8489]{Gerolymos_Senechal_Vallet_2009a}.
Recently, Arandiga \etal~\cite{Arandiga_Baeza_Belda_Mulet_2011a} gave analytical expressions of the optimal weights for the $(K_\mathrm{s}=r-1)$-level subdivision of
the upwind-biased stencil $\{i-(r-1),\cdots,i+(r-1)\}$.
Liu \etal~\cite{Liu_Shu_Zhang_2009a} showed that the weight-functions, in the reconstruction case, are rational functions,
expressed the weight-functions for the $(K_\mathrm{s}=\left\lceil\frac{M}{2}\right\rceil)$-level subdivision of the stencil $\{i-\lfloor\frac{M}{2}\rfloor,\cdots,i+M-\lfloor\frac{M}{2}\rfloor\}$,
in the range $M\in\{2,\cdots,11\}$, and computed the interval of convexity around $x_{i+\frac{1}{2}}$.
In \cite{Gerolymos_2012a}, we use the recurrence of \lemrefnp{Lem_GRRWFsMNARsAWENOID_s_GRRWFs_001} to study in detail the weight-functions for the Lagrange reconstructing polynomial \cite{Gerolymos_2011a},
obtain explicit recursive expressions for the weight-functions for an arbitrarily biased stencil on a homogeneous grid and an arbitrary subdivision level,
and determine the interval of convexity in the neighbourhood of $x_{i+\frac{1}{2}}$.
\qed
\end{remark}
%

Liu \etal~\cite{Liu_Shu_Zhang_2009a}, also study the representation of the first two derivatives ($n\in\{1,2\}$) of the Lagrange interpolating polynomial for the particular homogeneous stencils and subdivisions
studied in the reconstruction case \rmkref{Rmk_GRRWFsMNARsAWENOID_s_AnDLIP_001}.
In the present work, we show that the $(K_\mathrm{s}=1)$-level subdivision \defref{Def_GRRWFsMNARsAWENOID_s_I_002} weight-functions can be explicitly determined for the $n$-derivative of the
Lagrange interpolating polynomial \eqref{Eq_GRRWFsMNARsAWENOID_s_AnDLIP_001}. Using \lemrefnp{Lem_GRRWFsMNARsAWENOID_s_GRRWFs_001}, we define the weight-functions for the representation of the $n$-derivative of the
Lagrange interpolating polynomial \eqref{Eq_GRRWFsMNARsAWENOID_s_AnDLIP_001} on $\tsc{x}_{i-M_-,i+M_+}$ \defref{Def_GRRWFsMNARsAWENOID_s_I_001},
by combination of the $n$-derivative of the Lagrange interpolating polynomials on the $K_\mathrm{s}$-level substencils \defref{Def_GRRWFsMNARsAWENOID_s_I_002},
requiring that $n\leq M-K_\mathrm{s}$ so that the $n$-derivative be $\neq0_{\mathbb{R}_{M-K_\mathrm{s}-n}[x]}(x)$ on the substencils. The result is formulated in the following Proposition.
%
\begin{proposition}[{\rm Weight-functions for the $n$-derivative of the Lagrange interpolating polynomial}]
\label{Prp_GRRWFsMNARsAWENOID_s_AnDLIP_001}
\begin{subequations}
                                                                                                                                    \label{Eq_Prp_GRRWFsMNARsAWENOID_s_AnDLIP_001_001}
Assume the conditions of \prprefnp{Prp_GRRWFsMNARsAWENOID_s_ALIP_001}. Then, $\forall M_\pm\in{\mathbb Z}:\;M:=M_-+M_+\geq2$, $\forall K_\mathrm{s}\leq M-1$, $\forall n\leq M-K_\mathrm{s}$,
the $n$-derivative with respect to $x$ of the Lagrange interpolating polynomial on $\tsc{x}_{i-M_-,i+M_+}$ \defref{Def_GRRWFsMNARsAWENOID_s_I_001}
can be represented, almost everywhere, by combination of the $n$-derivative of the Lagrange interpolating polynomials on the $K_\mathrm{s}$-level substencils \defref{Def_GRRWFsMNARsAWENOID_s_I_002} of $\tsc{x}_{i-M_-,i+M_+}$, as
\begin{alignat}{6}
  p^{(n)}_{I,\tsc{x}_{i-M_-,i+M_+}}(x;f)= \sum_{k_\mathrm{s}=0}^{K_\mathrm{s}}\sigma_{I_{(n)},\tsc{x}_{i-M_-,i+M_+},K_\mathrm{s},k_\mathrm{s}}(x)
                                                                              \;p^{(n)}_{I,\tsc{x}_{i-M_-+k_\mathrm{s},i+M_+-K_\mathrm{s}+k_\mathrm{s}}}(x;f)
                                                                                    &\qquad\forall x\in{\mathbb R}\setminus\mathcal{S}_{I_{(n)},\tsc{x}_{i-M_-,i+M_+},K_\mathrm{s}}\quad
                                                                                                                                    \notag\\
                                                                                    &\qquad\forall f:{\mathbb R}\longrightarrow{\mathbb R}
                                                                                                                                    \label{Eq_Prp_GRRWFsMNARsAWENOID_s_AnDLIP_001_001a}
\end{alignat}
where the rational weight-functions $\sigma_{I_{(n)},\tsc{x}_{i-M_-,i+M_+},K_\mathrm{s},k_\mathrm{s}}(x)$ are defined recursively by
\begin{alignat}{6}
&\sigma_{I_{(n)},\tsc{x}_{i-M_-,i+M_+},K_\mathrm{s},k_\mathrm{s}}(x):=\left\{\begin{array}{ll}\dfrac{\alpha^{(n)}_{I,\tsc{x}_{i-M_-             ,i+M_+               },i-M_-+k_\mathrm{s}M}(x)}
                                                                                                    {\alpha^{(n)}_{I,\tsc{x}_{i-M_-+k_\mathrm{s},i+M_+-1+k_\mathrm{s}},i-M_-+k_\mathrm{s}M}(x)}      &K_\mathrm{s}=1       \\
                                                                                                                                                                    &                     \\
                                                                              {\displaystyle\sum_{\ell_\mathrm{s}=\max(0,k_\mathrm{s}-1)}^{\min(K_\mathrm{s}-1,k_\mathrm{s})}
                                                                                            \sigma_{I_{(n)},\tsc{x}_{i-M_-             ,i+M_+                                 },K_\mathrm{s}-1,             \ell_\mathrm{s}}(x)\;
                                                                                            \sigma_{I_{(n)},\tsc{x}_{i-M_-+\ell_{\rm s},i+M_+-(K_\mathrm{s}-1)+\ell_\mathrm{s}},             1,k_\mathrm{s}-\ell_\mathrm{s}}(x)
                                                                              }                                                                                     &K_\mathrm{s}\geq2    \\\end{array}\right.
                                                                                                                                    \notag\\
&\forall k_\mathrm{s}\in\{0,\cdots,K_\mathrm{s}\}\qquad\forall K_\mathrm{s}\in\{1,\cdots,M-1\}\qquad\forall n\in\{0,\cdots,M-K_\mathrm{s}\}
                                                                                                                                    \label{Eq_Prp_GRRWFsMNARsAWENOID_s_AnDLIP_001_001b}
\end{alignat}
and satisfy the consistency condition
\begin{alignat}{6}
 \sum_{k_\mathrm{s}=0}^{K_\mathrm{s}}\sigma_{I_{(n)},\tsc{x}_{i-M_-,i+M_+},K_\mathrm{s},k_\mathrm{s}}(x)=1\qquad\forall x\in{\mathbb R}\setminus\mathcal{S}_{I_{(n)},\tsc{x}_{i-M_-,i+M_+},K_\mathrm{s}}
                                                                                                                                    \label{Eq_Prp_GRRWFsMNARsAWENOID_s_AnDLIP_001_001c}
\end{alignat}
The set of poles of the rational weight-functions $\mathcal{S}_{I_{(n)},\tsc{x}_{i-M_-,i+M_+},K_\mathrm{s}}$ \eqref{Eq_Prp_GRRWFsMNARsAWENOID_s_AnDLIP_001_001a} satisfies
\begin{alignat}{6}
\mathcal{S}_{I_{(n)},\tsc{x}_{i-M_-,i+M_+},1           }\subseteq&\left\{x\in{\mathbb R}:\alpha^{(n)}_{I,\tsc{x}_{i-M_-+1,i+M_+},i+M_+}(x)=0\right\}
                                                                                                                                    \label{Eq_Prp_GRRWFsMNARsAWENOID_s_AnDLIP_001_001d}\\
\mathcal{S}_{I_{(n)},\tsc{x}_{i-M_-,i+M_+},K_\mathrm{s}}\subseteq&\bigcup_{   L_\mathrm{s}=0}^{K_\mathrm{s}-1}\displaylimits
                                                                  \bigcup_{\ell_\mathrm{s}=0}^{L_\mathrm{s}  }{\mathcal S}_{I_{(n)},\tsc{x}_{i-M_-+\ell_\mathrm{s},i+M_+-L_\mathrm{s}+\ell_\mathrm{s}},1}
                                                   =\left\{x\in{\mathbb R}:\prod_{   L_\mathrm{s}=0}^{K_\mathrm{s}-1}
                                                                           \prod_{\ell_\mathrm{s}=0}^{L_\mathrm{s}  }
                                                                           \alpha^{(n)}_{I,\tsc{x}_{i-M_-+1+\ell_\mathrm{s},i+M_+-L_\mathrm{s}+\ell_\mathrm{s}},i+M_+-L_\mathrm{s}+\ell_\mathrm{s}}(x)=0\right\}
                                                                                                                                    \notag\\
&\forall K_\mathrm{s}\in\{1,\cdots,M-1\}
                                                                                                                                    \label{Eq_Prp_GRRWFsMNARsAWENOID_s_AnDLIP_001_001e}
\end{alignat}
\end{subequations}
\end{proposition}
%
%
\begin{proof}
Because of \lemrefnp{Lem_GRRWFsMNARsAWENOID_s_GRRWFs_001} it suffices to prove that \prprefnp{Prp_GRRWFsMNARsAWENOID_s_AnDLIP_001} is valid for $K_\mathrm{s}=1$. Consider first the consistency
relation \eqref{Eq_Prp_GRRWFsMNARsAWENOID_s_AnDLIP_001_001c} for the $(K_\mathrm{s}=1)$-level weight-functions defined by \eqref{Eq_Prp_GRRWFsMNARsAWENOID_s_AnDLIP_001_001b},
\begin{subequations}
                                                                                                                                    \label{Eq_Prp_GRRWFsMNARsAWENOID_s_AnDLIP_001_002}
\begin{alignat}{6}
\sigma_{I_{(n)},\tsc{x}_{i-M_-,i+M_+},1,0}(x)=\dfrac{\alpha^{(n)}_{I,\tsc{x}_{i-M_-  ,i+M_+  },i-M_-}(x)}
                                                    {\alpha^{(n)}_{I,\tsc{x}_{i-M_-  ,i+M_+-1},i-M_-}(x)}
\qquad\forall n\in\{0,\cdots,M-1\}
                                                                                                                                    \label{Eq_Prp_GRRWFsMNARsAWENOID_s_AnDLIP_001_002a}\\
\sigma_{I_{(n)},\tsc{x}_{i-M_-,i+M_+},1,1}(x)=\dfrac{\alpha^{(n)}_{I,\tsc{x}_{i-M_-  ,i+M_+  },i+M_+}(x)}
                                                    {\alpha^{(n)}_{I,\tsc{x}_{i-M_-+1,i+M_+  },i+M_+}(x)}
\qquad\forall n\in\{0,\cdots,M-1\}
                                                                                                                                    \label{Eq_Prp_GRRWFsMNARsAWENOID_s_AnDLIP_001_002b}
\end{alignat}
\end{subequations}
By straightforward computation using the expression \eqref{Eq_Rmk_GRRWFsMNARsAWENOID_s_GRRWFs_003_001b} of the fundamental functions of Lagrange interpolation, we have
\begin{subequations}
                                                                                                                                    \label{Eq_Prp_GRRWFsMNARsAWENOID_s_AnDLIP_001_003}
\begin{alignat}{6}
\alpha_{I,\tsc{x}_{i-M_-  ,i+M_+-1},i-M_-}(x)= \prod_{k=-M_-+1}^{M_+-1}{\dfrac{x        -x_{i+k}}
                                                                              {x_{i-M_-}-x_{i+k}}}
                                             =&\left(\prod_{k=-M_-+1}^{M_+-1}{\dfrac{x_{i+M_+}-x_{i+k}}
                                                                                    {x_{i-M_-}-x_{i+k}}}\right)\prod_{k=-M_-+1}^{M_+-1}{\dfrac{x        -x_{i+k}}
                                                                                                                                              {x_{i+M_+}-x_{i+k}}}
                                                                                                                                    \notag\\
                                             =&\left(\prod_{k=-M_-+1}^{M_+-1}{\dfrac{x_{i+M_+}-x_{i+k}}
                                                                                    {x_{i-M_-}-x_{i+k}}}\right)\alpha_{I,\tsc{x}_{i-M_-+1,i+M_+  },i+M_+}(x)
                                                                                                                                    \label{Eq_Prp_GRRWFsMNARsAWENOID_s_AnDLIP_001_003a}
\end{alignat}
and taking into account that the product in the last line of \eqref{Eq_Prp_GRRWFsMNARsAWENOID_s_AnDLIP_001_003a} is independent of $x$, we have by differentiation
\begin{alignat}{6}
\alpha^{(n)}_{I,\tsc{x}_{i-M_-  ,i+M_+-1},i-M_-}(x)=\left(\prod_{k=-M_-+1}^{M_+-1}{\dfrac{x_{i+M_+}-x_{i+k}}
                                                                                         {x_{i-M_-}-x_{i+k}}}\right)\alpha^{(n)}_{I,\tsc{x}_{i-M_-+1,i+M_+  },i+M_+}(x)\quad\forall n\in\{0,\cdots,M-1\}
                                                                                                                                    \label{Eq_Prp_GRRWFsMNARsAWENOID_s_AnDLIP_001_003b}
\end{alignat}
where $n=0$ implies no differentiation, {\em ie} \eqref{Eq_Prp_GRRWFsMNARsAWENOID_s_AnDLIP_001_003a}. Since the product in \eqref{Eq_Prp_GRRWFsMNARsAWENOID_s_AnDLIP_001_003b} is independent of $n$, we also have
\begin{alignat}{6}
\dfrac{\alpha^{(n-1)}_{I,\tsc{x}_{i-M_-  ,i+M_+-1},i-M_-}(x)}
      {\alpha^{(n)  }_{I,\tsc{x}_{i-M_-  ,i+M_+-1},i-M_-}(x)}\stackrel{\eqref{Eq_Prp_GRRWFsMNARsAWENOID_s_AnDLIP_001_003b}}{=}
\dfrac{\alpha^{(n-1)}_{I,\tsc{x}_{i-M_-+1,i+M_+  },i+M_+}(x)}
      {\alpha^{(n)  }_{I,\tsc{x}_{i-M_-+1,i+M_+  },i+M_+}(x)}\quad\forall n\in\{1,\cdots,M-1\}
                                                                                                                                    \label{Eq_Prp_GRRWFsMNARsAWENOID_s_AnDLIP_001_003c}
\end{alignat}
We know by Aitken's Lemma \cite[pp. 204--205]{Henrici_1964a} that the $K_\mathrm{s}=1$ consistency relations hold for $n=0$, and by direct computation 
from the expression \eqref{Eq_Rmk_GRRWFsMNARsAWENOID_s_GRRWFs_003_001b} of the fundamental functions of Lagrange interpolation, we have \eqrefsab{Eq_Prp_GRRWFsMNARsAWENOID_s_ALIP_001_002a}
                                                                                                                                                 {Eq_Rmk_GRRWFsMNARsAWENOID_s_GRRWFs_003_002}
\begin{alignat}{6}
\dfrac{\alpha_{I,\tsc{x}_{i-M_-  ,i+M_+  },i-M_-}(x)}
      {\alpha_{I,\tsc{x}_{i-M_-  ,i+M_+-1},i-M_-}(x)}\stackrel{\eqrefsab{Eq_Prp_GRRWFsMNARsAWENOID_s_ALIP_001_002a}
                                                                        {Eq_Rmk_GRRWFsMNARsAWENOID_s_GRRWFs_003_002}}{=}&-\dfrac{x        -x_{i+M_+}}
                                                                                                                                {x_{i+M_+}-x_{i-M_-}}
                                                                                                                                    \label{Eq_Prp_GRRWFsMNARsAWENOID_s_AnDLIP_001_003d}\\
\dfrac{\alpha_{I,\tsc{x}_{i-M_-  ,i+M_+  },i+M_+}(x)}
      {\alpha_{I,\tsc{x}_{i-M_-+1,i+M_+  },i+M_+}(x)}\stackrel{\eqrefsab{Eq_Prp_GRRWFsMNARsAWENOID_s_ALIP_001_002a}
                                                                        {Eq_Rmk_GRRWFsMNARsAWENOID_s_GRRWFs_003_002}}{=}&+\dfrac{x        -x_{i-M_-}}
                                                                                                                                {x_{i+M_+}-x_{i-M_-}}
                                                                                                                                    \label{Eq_Prp_GRRWFsMNARsAWENOID_s_AnDLIP_001_003e}
\end{alignat}
which give by successive differentiation (proof by induction)
\begin{alignat}{6}
\alpha^{(n)}_{I,\tsc{x}_{i-M_-  ,i+M_+  },i-M_-}(x)\stackrel{\eqref{Eq_Prp_GRRWFsMNARsAWENOID_s_AnDLIP_001_003d}}{=}&-\dfrac{x        -x_{i+M_+}}
                                                                                                                            {x_{i+M_+}-x_{i-M_-}}\alpha^{(n)  }_{I,\tsc{x}_{i-M_-  ,i+M_+-1},i-M_-}(x)
                                                                                                                     -\dfrac{n                  }
                                                                                                                            {x_{i+M_+}-x_{i-M_-}}\alpha^{(n-1)}_{I,\tsc{x}_{i-M_-  ,i+M_+-1},i-M_-}(x)
                                                                                                                                    \label{Eq_Prp_GRRWFsMNARsAWENOID_s_AnDLIP_001_003f}\\
\alpha^{(n)}_{I,\tsc{x}_{i-M_-  ,i+M_+  },i+M_+}(x)\stackrel{\eqref{Eq_Prp_GRRWFsMNARsAWENOID_s_AnDLIP_001_003d}}{=}&+\dfrac{x        -x_{i-M_-}}
                                                                                                                            {x_{i+M_+}-x_{i-M_-}}\alpha^{(n)  }_{I,\tsc{x}_{i-M_-+1,i+M_+  },i+M_+}(x)
                                                                                                                     +\dfrac{n                  }
                                                                                                                            {x_{i+M_+}-x_{i-M_-}}\alpha^{(n-1)}_{I,\tsc{x}_{i-M_-+1,i+M_+  },i+M_+}(x)
                                                                                                                                    \label{Eq_Prp_GRRWFsMNARsAWENOID_s_AnDLIP_001_003g}
\end{alignat}
\end{subequations}
Using the above relations \eqrefsab{Eq_Prp_GRRWFsMNARsAWENOID_s_AnDLIP_001_002}{Eq_Prp_GRRWFsMNARsAWENOID_s_AnDLIP_001_003}, we readily have
\begin{alignat}{6}
\sigma_{I_{(n)},\tsc{x}_{i-M_-,i+M_+},1,0}(x)+
\sigma_{I_{(n)},\tsc{x}_{i-M_-,i+M_+},1,1}(x)\stackrel{\eqrefsabc{Eq_Prp_GRRWFsMNARsAWENOID_s_AnDLIP_001_002}
                                                                 {Eq_Prp_GRRWFsMNARsAWENOID_s_AnDLIP_001_003f}
                                                                 {Eq_Prp_GRRWFsMNARsAWENOID_s_AnDLIP_001_003g}}{=}&\underbrace{-\dfrac{x        -x_{i+M_+}}
                                                                                                                                      {x_{i+M_+}-x_{i-M_-}}
                                                                                                                               +\dfrac{x        -x_{i-M_-}}
                                                                                                                                      {x_{i+M_+}-x_{i-M_-}}}_{=1}
                                                                                                                                    \notag\\
                                                                                                                  &-\dfrac{n                  }
                                                                                                                          {x_{i+M_+}-x_{i-M_-}}\underbrace{\left(\dfrac{\alpha^{(n-1)}_{I,\tsc{x}_{i-M_-  ,i+M_+-1},i-M_-}(x)}
                                                                                                                                                                       {\alpha^{(n)  }_{I,\tsc{x}_{i-M_-  ,i+M_+-1},i-M_-}(x)}
                                                                                                                                                                -\dfrac{\alpha^{(n-1)}_{I,\tsc{x}_{i-M_-+1,i+M_+  },i+M_+}(x)}
                                                                                                                                                                       {\alpha^{(n)  }_{I,\tsc{x}_{i-M_-+1,i+M_+  },i+M_+}(x)}\right)
                                                                                                                                                          }_{\stackrel{\eqref{Eq_Prp_GRRWFsMNARsAWENOID_s_AnDLIP_001_003c}}{=}0}
                                                                                                                                    \notag\\
=&1\qquad\forall n\in\{1,\cdots,M-1\}
                                                                                                                                    \label{Eq_Prp_GRRWFsMNARsAWENOID_s_AnDLIP_001_004}
\end{alignat}
By \eqref{Eq_Prp_GRRWFsMNARsAWENOID_s_AnDLIP_001_004}, we have proven the consistency condition \eqref{Eq_Prp_GRRWFsMNARsAWENOID_s_AnDLIP_001_001c}, for $K_\mathrm{s}=1$, $\forall n\in\{0,\cdots,M-1\}$,
the case $n=0$ following from \prprefnp{Prp_GRRWFsMNARsAWENOID_s_ALIP_001}. Obviously the weight-functions \eqref{Eq_Prp_GRRWFsMNARsAWENOID_s_AnDLIP_001_002} are defined almost everywhere,
except at the roots of the denominator, which because of \eqref{Eq_Prp_GRRWFsMNARsAWENOID_s_AnDLIP_001_003b}, are defined by \eqref{Eq_Prp_GRRWFsMNARsAWENOID_s_AnDLIP_001_001d}. The $\subseteq$ relation
in \eqref{Eq_Prp_GRRWFsMNARsAWENOID_s_AnDLIP_001_001d} is used, because there may be cancelation of poles by polynomial division, as in the interpolation case ($n=0$),
where there are no singularities \eqref{Eq_Prp_GRRWFsMNARsAWENOID_s_ALIP_001_001}.

To prove the representation \eqref{Eq_Prp_GRRWFsMNARsAWENOID_s_AnDLIP_001_001a}, for $K_\mathrm{s}=1$, we start from the remainder theorem of the Lagrange interpolating polynomial,
which states \cite[Theorem 9.2, p. 187]{Henrici_1964a} that for any real function $f:\mathbb{R}\to\mathbb{R}$ of class $C^{M+1}$
\begin{subequations}
                                                                                                                                    \label{Eq_Prp_GRRWFsMNARsAWENOID_s_AnDLIP_001_005}
\begin{alignat}{6}
\forall f\in C^{M+1}(\mathbb{R})\quad\forall x\in\mathbb{R}\quad\exists\;t(x;\tsc{x}_{i-M_-,i+M_+};f)\in[\min(x,x_{i-M_-}),\max(x,x_{i+M_+})]\;:\;
                                                                                                                                    \notag\\
p_{I,\tsc{x}_{i-M_-,i+M_+}}(x;f)=f(x)-\dfrac{1}{(M+1)!}\left(\prod_{k=-M_-}^{M_+}(x-x_{i+k})\right)f^{(M+1)}(t)
                                                                                                                                    \label{Eq_Prp_GRRWFsMNARsAWENOID_s_AnDLIP_001_005a}
\end{alignat}
Since for any polynomial of degree $\leq M$ we have
\begin{alignat}{6}
\forall q\in\mathbb{R}_M[x]\Longrightarrow\left\{\begin{array}{l}q^{(M)}(x)=\mathrm{coeff}[x^M,q(x)]\;M!\\
                                                                 q^{(M+1)}(x)=0                         \\\end{array}\right.\qquad\forall x\in\mathbb{R}
                                                                                                                                    \label{Eq_Prp_GRRWFsMNARsAWENOID_s_AnDLIP_001_005b}
\end{alignat}
we can write because of \eqref{Eq_Prp_GRRWFsMNARsAWENOID_s_AnDLIP_001_005a}
\begin{alignat}{6}
p_{I,\tsc{x}_{i-M_-  ,i+M_+-1}}(x;q)=&q(x)-\prod_{k=-M_-}^{M_+-1}(x-x_{i+k})
                                                                                                                                    \label{Eq_Prp_GRRWFsMNARsAWENOID_s_AnDLIP_001_005c}\\
p_{I,\tsc{x}_{i-M_-  ,i+M_+  }}(x;q)=&q(x)                                               &\qquad\forall q\in\mathbb{R}_M[x]\;:\;\mathrm{coeff}[x^M,q(x)]=1
                                                                                                                                    \label{Eq_Prp_GRRWFsMNARsAWENOID_s_AnDLIP_001_005d}\\
p_{I,\tsc{x}_{i-M_-+1,i+M_+  }}(x;q)=&q(x)-\prod_{k=-M_-+1}^{M_+}(x-x_{i+k})
                                                                                                                                    \label{Eq_Prp_GRRWFsMNARsAWENOID_s_AnDLIP_001_005e}
\end{alignat}
Using the expression \eqref{Eq_Rmk_GRRWFsMNARsAWENOID_s_GRRWFs_003_001b} of the fundamental functions of Lagrange interpolation
in \eqrefsab{Eq_Prp_GRRWFsMNARsAWENOID_s_AnDLIP_001_005c}
            {Eq_Prp_GRRWFsMNARsAWENOID_s_AnDLIP_001_005e},\footnote{\label{fn_Prp_GRRWFsMNARsAWENOID_s_AnDLIP_001_001}\begin{alignat}{6}
                                                                                                                      \prod_{k=-M_-}^{M_+-1}(x-x_{i+k})\stackrel{\eqref{Eq_Rmk_GRRWFsMNARsAWENOID_s_GRRWFs_003_001b}}{=}
                                                                                                                      \left(\prod_{k=-M_-}^{M_+-1}(x_{i+M_+}-x_{i+k})\right)\alpha_{I,\tsc{x}_{i-M_-  ,i+M_+  },i+M_+}(x)
                                                                                                                      \notag\\
                                                                                                                      \prod_{k=-M_-+1}^{M_+}(x-x_{i+k})\stackrel{\eqref{Eq_Rmk_GRRWFsMNARsAWENOID_s_GRRWFs_003_001b}}{=}
                                                                                                                      \left(\prod_{k=-M_-+1}^{M_+}(x_{i-M_-}-x_{i+k})\right)\alpha_{I,\tsc{x}_{i-M_-  ,i+M_+  },i-M_-}(x)
                                                                                                                      \notag
                                                                                                                      \end{alignat}
                                                                   }
and differentiating, we have, $\forall n\in\mathbb{N}_0$,
\begin{alignat}{6}
p^{(n)}_{I,\tsc{x}_{i-M_-  ,i+M_+-1}}(x;q)=&q^{(n)}(x)-\left(\prod_{k=-M_-}^{M_+-1}(x_{i+M_+}-x_{i+k})\right)\alpha^{(n)}_{I,\tsc{x}_{i-M_-  ,i+M_+  },i+M_+}(x)
                                                                                                                                    \label{Eq_Prp_GRRWFsMNARsAWENOID_s_AnDLIP_001_005f}\\
p^{(n)}_{I,\tsc{x}_{i-M_-  ,i+M_+  }}(x;q)=&q^{(n)}(x)                                               &\quad\forall q\in\mathbb{R}_M[x]\;:\;\mathrm{coeff}[x^M,q(x)]=1
                                                                                                                                    \label{Eq_Prp_GRRWFsMNARsAWENOID_s_AnDLIP_001_005g}\\
p^{(n)}_{I,\tsc{x}_{i-M_-+1,i+M_+  }}(x;q)=&q^{(n)}(x)-\left(\prod_{k=-M_-+1}^{M_+}(x_{i-M_-}-x_{i+k})\right)\alpha^{(n)}_{I,\tsc{x}_{i-M_-  ,i+M_+  },i-M_-}(x)
                                                                                                                                    \label{Eq_Prp_GRRWFsMNARsAWENOID_s_AnDLIP_001_005h}
\end{alignat}
Combining \eqrefsabc{Eq_Prp_GRRWFsMNARsAWENOID_s_AnDLIP_001_002}
                    {Eq_Prp_GRRWFsMNARsAWENOID_s_AnDLIP_001_005f}
                    {Eq_Prp_GRRWFsMNARsAWENOID_s_AnDLIP_001_005h}, and using \eqref{Eq_Prp_GRRWFsMNARsAWENOID_s_AnDLIP_001_003b}, we have
\begin{alignat}{6}
&\sigma_{I_{(n)},\tsc{x}_{i-M_-,i+M_+},1,0}(x)\;p^{(n)}_{I,\tsc{x}_{i-M_-  ,i+M_+-1}}(x;q)+
 \sigma_{I_{(n)},\tsc{x}_{i-M_-,i+M_+},1,1}(x)\;p^{(n)}_{I,\tsc{x}_{i-M_-+1,i+M_+  }}(x;q)
                                                                                                                                    \notag\\
\stackrel{\eqrefsabc{Eq_Prp_GRRWFsMNARsAWENOID_s_AnDLIP_001_002}
                    {Eq_Prp_GRRWFsMNARsAWENOID_s_AnDLIP_001_005f}
                    {Eq_Prp_GRRWFsMNARsAWENOID_s_AnDLIP_001_005h}}{=}
 &\underbrace{\left(\sigma_{I_{(n)},\tsc{x}_{i-M_-,i+M_+},1,0}(x)+\sigma_{I_{(n)},\tsc{x}_{i-M_-,i+M_+},1,1}(x)\right)}_{\stackrel{\eqref{Eq_Prp_GRRWFsMNARsAWENOID_s_AnDLIP_001_004}}{=}1}q^{(n)}(x)
                                                                                                                                    \notag\\
-&(x_{i+M_+}-x_{i-M_-})\;\alpha^{(n)}_{I,\tsc{x}_{i-M_-  ,i+M_+  },i+M_+}(x)\;\alpha^{(n)}_{I,\tsc{x}_{i-M_-  ,i+M_+  },i-M_-}(x)
\underbrace{\left(\dfrac{\displaystyle\prod_{k=-M_-+1}^{M_+-1}(x_{i+M_+}-x_{i+k})}
                        {\alpha^{(n)}_{I,\tsc{x}_{i-M_-  ,i+M_+-1},i-M_-}(x)     }
                 -\dfrac{\displaystyle\prod_{k=-M_-+1}^{M_+-1}(x_{i-M_-}-x_{i+k})}
                        {\alpha^{(n)}_{I,\tsc{x}_{i-M_-+1,i+M_+  },i+M_+}(x)     }\right)}_{\stackrel{\eqref{Eq_Prp_GRRWFsMNARsAWENOID_s_AnDLIP_001_003b}}{=}0}
                                                                                                                                    \notag\\
=&q^{(n)}(x)\stackrel{\eqref{Eq_Prp_GRRWFsMNARsAWENOID_s_AnDLIP_001_005g}}{=}p^{(n)}_{I,\tsc{x}_{i-M_-  ,i+M_+  }}(x;q)\quad\forall q\in\mathbb{R}_M[x]\;:\;\mathrm{coeff}[x^M,q(x)]=1\quad\forall n\in\{1,\cdots,M-1\}
                                                                                                                                    \label{Eq_Prp_GRRWFsMNARsAWENOID_s_AnDLIP_001_005i}
\end{alignat}
Using the representation \eqref{Eq_GRRWFsMNARsAWENOID_s_AnDLIP_001} of the $n$-derivative of Langrange interpolating polynomials in \eqref{Eq_Prp_GRRWFsMNARsAWENOID_s_AnDLIP_001_005i}, we readily have
\begin{alignat}{6}
 &\sigma_{I_{(n)},\tsc{x}_{i-M_-,i+M_+},1,0}(x)\;\alpha^{(n)}_{I,\tsc{x}_{i-M_-  ,i+M_+-1},i-M_-}(x)\;q(x_{i-M_-})
                                                                                                                                    \notag\\
+&\sum_{\ell=-M_-+1}^{M_+-1}\Big(\sigma_{I_{(n)},\tsc{x}_{i-M_-,i+M_+},1,0}(x)\;\alpha^{(n)}_{I,\tsc{x}_{i-M_-  ,i+M_+-1},i+\ell}(x)
                                +\sigma_{I_{(n)},\tsc{x}_{i-M_-,i+M_+},1,1}(x)\;\alpha^{(n)}_{I,\tsc{x}_{i-M_-+1,i+M_+  },i+\ell}(x)\Big)\;q(x_{i+\ell})
                                                                                                                                    \notag\\
+&\sigma_{I_{(n)},\tsc{x}_{i-M_-,i+M_+},1,1}(x)\;\alpha^{(n)}_{I,\tsc{x}_{i-M_-+1,i+M_+  },i+M_+}(x)\;q(x_{i+M_+})
                                                                                                                                    \notag\\
=&\sum_{\ell=-M_-}^{M_+}\alpha^{(n)}_{I,\tsc{x}_{i-M_-,i+M_+},i+\ell}(x)\;q(x_{i+\ell})\qquad\forall q\in\mathbb{R}_M[x]\;:\;\mathrm{coeff}[x^M,q(x)]=1\quad\forall n\in\{1,\cdots,M-1\}
                                                                                                                                    \label{Eq_Prp_GRRWFsMNARsAWENOID_s_AnDLIP_001_005j}
\end{alignat}
Applying \eqref{Eq_Prp_GRRWFsMNARsAWENOID_s_AnDLIP_001_005j}, successively, to the polynomials
\begin{alignat}{6}
\mathbb{R}_M[x]\ni
\prod_{\substack{m=-M_-\\
                 m\neq k}  }^{M_+}(x-x_{i+m})=0\quad\forall x\in\left\{x_{i-M_-},\cdots,x_{i+M_+}\right\}\setminus\{x_{i+k}\}\quad\forall k\in\{-M_-,\cdots,+M_+\}
                                                                                                        \label{Eq_Lem_RLRPCS_s_RCsSs_ss_Ks1_001_003i}
\end{alignat}
yields
\begin{alignat}{6}
\alpha^{(n)}_{I,\tsc{x}_{i-M_-,i+M_+},i-M_- }(x)=&\sigma_{I_{(n)},\tsc{x}_{i-M_-,i+M_+},1,0}(x)\;\alpha^{(n)}_{I,\tsc{x}_{i-M_-  ,i+M_+-1},i-M_-}(x)\quad\iff\eqref{Eq_Prp_GRRWFsMNARsAWENOID_s_AnDLIP_001_002a}
                                                                                                                                    \label{Eq_Prp_GRRWFsMNARsAWENOID_s_AnDLIP_001_005k}\\
\alpha^{(n)}_{I,\tsc{x}_{i-M_-,i+M_+},i+\ell}(x)=&\sigma_{I_{(n)},\tsc{x}_{i-M_-,i+M_+},1,0}(x)\;\alpha^{(n)}_{I,\tsc{x}_{i-M_-  ,i+M_+-1},i+\ell}(x)
                                                                                                                                    \notag\\
                                                +&\sigma_{I_{(n)},\tsc{x}_{i-M_-,i+M_+},1,1}(x)\;\alpha^{(n)}_{I,\tsc{x}_{i-M_-+1,i+M_+  },i+\ell}(x)\qquad\forall\ell\in\{-M_-+1,\cdots,M_+-1\}
                                                                                                                                    \label{Eq_Prp_GRRWFsMNARsAWENOID_s_AnDLIP_001_005l}\\
\alpha^{(n)}_{I,\tsc{x}_{i-M_-,i+M_+},i+M_+ }(x)=&\sigma_{I_{(n)},\tsc{x}_{i-M_-,i+M_+},1,1}(x)\;\alpha^{(n)}_{I,\tsc{x}_{i-M_-+1,i+M_+  },i+M_+}(x)\quad\iff\eqref{Eq_Prp_GRRWFsMNARsAWENOID_s_AnDLIP_001_002b}
                                                                                                                                    \label{Eq_Prp_GRRWFsMNARsAWENOID_s_AnDLIP_001_005m}
\end{alignat}
\end{subequations}
Replacing the expression \eqref{Eq_GRRWFsMNARsAWENOID_s_AnDLIP_001} for the $n$-derivative of Langrange interpolating polynomials in \eqref{Eq_Prp_GRRWFsMNARsAWENOID_s_AnDLIP_001_001a},
proves, by \eqrefsatob{Eq_Prp_GRRWFsMNARsAWENOID_s_AnDLIP_001_005k}
                      {Eq_Prp_GRRWFsMNARsAWENOID_s_AnDLIP_001_005m},
\eqref{Eq_Prp_GRRWFsMNARsAWENOID_s_AnDLIP_001_001a} $\forall f:\mathbb{R}\to\mathbb{R}$, for $K_\mathrm{s}=1$.
\lemrefnp{Lem_GRRWFsMNARsAWENOID_s_GRRWFs_001} completes the proof.
\qed
\end{proof}
%
%
\begin{remark}[{\rm Relation of \prprefnp{Prp_GRRWFsMNARsAWENOID_s_AnDLIP_001} to previous work}]
\label{Rmk_GRRWFsMNARsAWENOID_s_AnDLIP_002}
To the author's knowledge, the case $n\geq3$ has not been studied before. For $n\in\{1,2\}$, the simple analytical recursive expression \eqref{Eq_Prp_GRRWFsMNARsAWENOID_s_AnDLIP_001_001b}
for the weight-functions agrees with the expressions obtained by Liu \etal~\cite{Liu_Shu_Zhang_2009a},
for the case of a homogeneous grid, using symbolic calculation for $r\in\{2,\cdots,7\}$,
for the ($K_\mathrm{s}=r-1$)-level subdivision of the stencils
$\tsc{x}_{i-r,i+(r-1)}:=\{-r+\tfrac{1}{2},\cdots,r-\tfrac{1}{2}\}$ \cite[Tab. 3.2, p. 516, for $n=1$, and Tab. 3.8, p. 520, for $n=2$]{Liu_Shu_Zhang_2009a} and
for the ($K_\mathrm{s}=r$)-level subdivision of the stencils
$\tsc{x}_{i-r,i+r}:=\{-r+\tfrac{1}{2},\cdots,r+\tfrac{1}{2}\}$ \cite[Tab. 3.5, p. 518, for $n=1$, and Tab. 3.9, p. 521, for $n=2$]{Liu_Shu_Zhang_2009a}.
\prprefnp{Prp_GRRWFsMNARsAWENOID_s_AnDLIP_001} formally proves the existence and analytical expression of the weight-functions, for the $n$-derivative ($\forall n\in\{0,\cdots,M-K_\mathrm{s}\}$),
for a general $K_\mathrm{s}$-level subdivision ($K_\mathrm{s}\in\{1,\cdots,M-1\}$), for an arbitrary stencil of $M+1$ distinct ordered points \defref{Def_GRRWFsMNARsAWENOID_s_I_001},
on an inhomogeneous grid.
\qed
\end{remark}
%

%
%
%
%
%
%
%
%
%
\section{Conclusions}\label{GRRWFsMNARsAWENOID_s_C}
%
%
%
%
%
%
%
%
%

Every system of functions \eqref{Eq_Lem_GRRWFsMNARsAWENOID_s_GRRWFs_001_001a}, depending on 2 integer parameters,
which is equipped with an associated system of weight-functions satisfying a $1$-level subdivision property \eqrefsab{Eq_Lem_GRRWFsMNARsAWENOID_s_GRRWFs_001_001b}
                                                                                                                     {Eq_Lem_GRRWFsMNARsAWENOID_s_GRRWFs_001_001c}
also satisfies $K_\mathrm{s}$-level subdivision relations \lemref{Lem_GRRWFsMNARsAWENOID_s_GRRWFs_001},
with weight-functions generated by the recurrence \eqref{Eq_Lem_GRRWFsMNARsAWENOID_s_GRRWFs_001_001e},
which can be interpreted as an inverted generalized Neville algorithm \cite[pp. 207--208]{Henrici_1964a}.

As an application of \lemrefnp{Lem_GRRWFsMNARsAWENOID_s_GRRWFs_001} we developed simple explicit expressions for $K_\mathrm{s}$-level weight-functions
of the Lagrange interpolating polynomial \prpref{Prp_GRRWFsMNARsAWENOID_s_ALIP_001} on a general stencil in an inhomogeneous grid, which allow explicit determination of the interval
of positivity of the weight-functions \prpref{Prp_GRRWFsMNARsAWENOID_s_ALIP_002} generalizing previous results \cite{Carlini_Ferretti_Russo_2005a,
                                                                                                                     Liu_Shu_Zhang_2009a}.
By \eqref{Eq_Prp_GRRWFsMNARsAWENOID_s_ALIP_002_001} the length of the positivity interval is $M_+-M_--2K_\mathrm{s}+2$ cells.

We further investigated the existence of $K_\mathrm{s}$-level weight-functions for the $n$-derivative of the Lagrange interpolating polynomial ($n\in\{1,\cdots,M-K_\mathrm{s}\}$).
Having proved simple analytical expressions for the $1$-level weight-functions,
\lemrefnp{Lem_GRRWFsMNARsAWENOID_s_GRRWFs_001} was applied to develop an analytical recursive expression for the $K_\mathrm{s}$-level weight-functions
of the $n$-derivative ($n\in\{1,\cdots,M-K_\mathrm{s}\}$). These results are valid for general inhomogeneous grids.

Other potential applications of \lemrefnp{Lem_GRRWFsMNARsAWENOID_s_GRRWFs_001} include \tsc{weno} integration \cite{Liu_Shu_Zhang_2009a},
the development of \tsc{weno} schemes for biased near-boundary stencils, and other than Lagrange types of polynomial interpolation \cite{Shu_2009a}.

%
%
%
%
%
%
%
%
%
\section*{Acknowledgments}
%
%
%
%
%
%
%
%
%

Computations were performed using \tsc{hpc} resources from \tsc{genci--idris} (Grants 2010--066327 and 2011--066327).

%
%
%
%
%
%
%
%
%
\footnotesize\bibliographystyle{elsarticle-num}\bibliography{Aerodynamics,GV_news}
%
%
%
%
%
%
%
%
%

\end{document}